# EXACT HAUSDORFF MEASURE ON THE BOUNDARY OF A GALTON–WATSON TREE

By Toshiro Watanabe

*University of Aizu*

A necessary and sufficient condition for the almost sure existence of an absolutely continuous (with respect to the branching measure) exact Hausdorff measure on the boundary of a Galton–Watson tree is obtained. In the case where the absolutely continuous exact Hausdorff measure does not exist almost surely, a criterion which classifies gauge functions $\phi$ according to whether $\phi$-Hausdorff measure of the boundary minus a certain exceptional set is zero or infinity is given. Important examples are discussed in four additional theorems. In particular, Hawkes's conjecture in 1981 is solved. Problems of determining the exact local dimension of the branching measure at a typical point of the boundary are also solved.

**1. Introduction.** An interesting history of the classical problem of determining the Hausdorff and packing dimensions and then the exact Hausdorff and packing measures of the boundary of a supercritical Galton–Watson tree is found in the previous paper [46]. It was initiated in 1973 by the thesis of Holmes [18], whose supervisor and examiner were C. A. Rogers and S. J. Taylor, respectively. The author [46] completely solved the problem of determining the exact packing measure of the boundary of the tree by filling the critical gap in the proof of the theorem of Liu [22], which had been pointed out by Berlinkov and Mauldin [4]. Berlinkov [3] independently studied the exact packing measures of homogeneous random recursive fractals and, as a corollary, he obtained an analogous result under a certain additional assumption on the tree. However, it was stated without precise proof and he could not identify the explicit value of the exact packing measure of the boundary. Upon an outline of Hawkes [17], the author [46] defined a random sequence $\{Y(n)\}$ for $n \leq 0$ as $Y(-n) := \mu(B_n)$, that is, the branching measure of the ball $B_n$ with diameter $e^{-n}$ on the boundary of the tree and discovered that it is a shift self-similar additive random sequence on a certain extended probability space. It is a key fact for solving this old problem, which enables us to use new limit theorems for shift self-similar additive random sequences developed by the author [44]. In the present paper, we extensively employ limit theorems of "limsup" type for the sequence $\{Y(n)\}$ and find a necessary and sufficient condition for the almost sure existence of an absolutely continuous (with respect to the branching measure) exact Hausdorff measure on the boundary of a Galton–Watson tree. It is represented by the nondominated variation of the right tail of the martingale limit of the branching process, equivalently, by the nondominated variation of the integrated function of the right tail of the offspring distribution. See Corollary 1.1. In the case where an absolutely continuous exact Hausdorff measure does not exist, we give a criterion which classifies gauge functions $\phi$ according to whether $\phi$-Hausdorff measure of the boundary minus a certain exceptional set is 0 or $\infty$. See Theorem 1.2. The explicit value of $\phi$-Hausdorff measure of the boundary is determined for each example in three additional theorems by closing the serious gaps in the proofs of Liu [21]. See Remark 1.4. In particular, Theorem 1.3 can be applied to obtain upper and lower bounds for the explicit value of the exact Hausdorff measure of a homogeneous random recursive fractal such as the limit set of Mandelbrot's fractal percolation and the path of a self-avoiding process on the Sierpinski gasket. See [14] and the examples of Berlinkov [3]. Moreover, a conjecture of Hawkes [17] in 1981 is solved. See Theorem 1.6. As is found in the concluding remarks, our problem of determining the exact Hausdorff measure is not yet completely solved. However, it is realized that the study of the exceptional set $\Delta$ defined by (1.7) below will lead to the complete solution.

In what follows, denote by $\mathbb{R}^d$ the $d$-dimensional Euclidean space and let $\mathbb{R}_+ = [0, \infty)$. Let $\mathbb{Z} = \{0, \pm 1, \pm 2, \ldots\}$, $\mathbb{Z}_+ = \{0, 1, 2, \ldots\}$, $\mathbb{N} = \{1, 2, 3, \ldots\}$, and denote $\mathbf{U} = \bigcup_{n=0}^{\infty} \mathbb{Z}_+^n$ with $\mathbb{Z}_+^0 = \varnothing$. We denote $\mathbf{i} \in \mathbb{Z}_+^n$ by $(i_k)_{k=1}^n$ or $(i_1, i_2, \ldots, i_n)$. For $\mathbf{i} \in \mathbb{Z}_+^n \subset \mathbf{U}$, we define $|\mathbf{i}| = n$. Let $\mathbf{I} = \mathbb{Z}_+^{\mathbb{N}}$. We denote $\mathbf{i} \in \mathbf{I}$ by $(i_k)_{k=1}^{\infty}$ and define, for $\mathbf{i} \in \mathbf{I}$, $|\mathbf{i}| = \infty$. For $\mathbf{i} = (i_k)_{k=1}^n$ and $\mathbf{j} = (j_k)_{k=1}^m$ in $\mathbf{U}$, we define $\mathbf{i} * \mathbf{j} \in \mathbf{U}$ as $\mathbf{i} * \mathbf{j} := (i_1, i_2, \ldots, i_n, j_1, \ldots, j_m)$. In particular, we have $\varnothing * \mathbf{i} = \mathbf{i} * \varnothing = \mathbf{i}$. We define $\mathbf{i}|n = (i_k)_{k=1}^n$ for $\mathbf{i} \in \mathbf{U} \cup \mathbf{I}$ with $n \in \mathbb{Z}_+ \cup \{\infty\}$ satisfying $n \leq |\mathbf{i}|$. We understand that $\mathbf{i}|0 = \varnothing$. We say that $\mathbf{i} \leq \mathbf{j}$ in $\mathbf{U} \cup \mathbf{I}$ if $|\mathbf{i}| = n \leq |\mathbf{j}|$ and $\mathbf{j}|n = \mathbf{i}$. In this order, we define $\mathbf{i} \wedge \mathbf{j} \in \mathbf{U} \cup \mathbf{I}$ for $\mathbf{i}, \mathbf{j} \in \mathbf{U} \cup \mathbf{I}$ as $\mathbf{i} \wedge \mathbf{j} := \max\{\mathbf{k} \in \mathbf{U} \cup \mathbf{I} : \mathbf{k} \leq \mathbf{i} \text{ and } \mathbf{k} \leq \mathbf{j}\}$. We define a metric $d(\mathbf{i}, \mathbf{j})$ for $\mathbf{i}, \mathbf{j} \in \mathbf{I}$ as $d(\mathbf{i}, \mathbf{j}) := e^{-|\mathbf{i} \wedge \mathbf{j}|}$. Then $(\mathbf{I}, d)$ is an ultrametric space. Denote by $\mathcal{B}(\mathbf{I})$ the class of all Borel sets in $(\mathbf{I}, d)$. From now on, let $\{N_{\mathbf{i}}, \mathbf{i} \in \mathbf{U}\}$ be $\mathbb{Z}_+$-valued i.i.d. random variables on a probability space $(\Omega, \mathcal{F}, P)$. In particular, put $N := N_{\varnothing}$. We assume, to avoid the trivial cases, that the support of the distribution of $N$ is not a one-point set. We denote by $f(s) := \sum_{n=0}^{\infty} p_n s^n$ the probability generating function (p.g.f. for short) of the distribution of $N$, where $p_n := P(N = n)$ for $n \in \mathbb{Z}_+$. Let $f_n(s)$ be the $n$th iteration of $f(s)$ with itself. We assume that

$$(1.1) \qquad a := E(N) > 1 \quad \text{and} \quad E(N \log N) < \infty.$$

A set $\mathbf{T} \subset \mathbf{U}$ is called a *Galton–Watson tree* on $(\Omega, \mathcal{F}, P)$ with *offspring distribution* $\{p_n\}_{n \geq 0} = \{P(N = n)\}_{n \geq 0}$ if the following three conditions are satisfied:

(1) $\varnothing \in \mathbf{T}$.
(2) Let $\mathbf{i} \in \mathbf{T}$ and $i \in \mathbb{Z}_+$. Then $\mathbf{i} * i \in \mathbf{T}$ if and only if $0 \leq i \leq N_{\mathbf{i}} - 1$.
(3) If $\mathbf{i} \in \mathbf{T}$ and $\mathbf{j} \leq \mathbf{i}$, then $\mathbf{j} \in \mathbf{T}$.

Let $\mathbf{T}$ be a Galton–Watson tree. We define the *boundary* (or *branching set*) $\partial \mathbf{T}$ of $\mathbf{T}$ as

$$\partial \mathbf{T} := \{\mathbf{i} \in \mathbf{I} : \mathbf{i} | n \in \mathbf{T} \text{ for every } n \in \mathbb{Z}_+\}.$$

We define $F_n \subset \mathbf{T}$ and $Z_n$ for $n \in \mathbb{Z}_+$ as

$$F_n := \{\mathbf{i} \in \mathbf{T} : |\mathbf{i}| = n\} \quad \text{and} \quad Z_n := \operatorname{Card} F_n.$$

Here $\operatorname{Card} A$ stands for the cardinality of a set $A$. Then $\{Z_n, n \in \mathbb{Z}_+\}$ is a supercritical Galton–Watson branching process with p.g.f. $f(s)$ of the number of offspring. Thus the distribution of $N$ is the same as that of $Z_1$. We define an $\mathbb{R}_+$-valued random variable $W$ as the following martingale limit:

$$(1.2) \qquad W := \lim_{n \to \infty} \frac{Z_n}{a^n}.$$

Our assumption (1.1) implies that $W$ exists almost surely with $E(W) = 1$. Note that $\{W = 0\} = \{\partial \mathbf{T} = \varnothing\} = \{\lim_{n \to \infty} Z_n = 0\}$ up to probability zero sets and that $q := P(W = 0)$ is the first nonnegative solution of the equation $f(s) - s = 0$. Thus it is obvious that $q = 0$ if and only if $p_0 = 0$. See [1] as to the above assertions. In this paper we use the words "increase" and "decrease" in the wide sense allowing flatness. A nonnegative decreasing function $h(x)$ on $\mathbb{R}_+$ is called of *dominated variation* as $x \to \infty$ [$h(x) \in \mathcal{D}$ for short] if $h(x) > 0$ for $x > 0$ and $\liminf_{x \to \infty} h(2x)/h(x) > 0$. Note that if $h(x)$ is regularly varying as $x \to \infty$, then $h(x) \in \mathcal{D}$. See [10] and [39]. For two positive functions $h_1(x)$ and $h_2(x)$ on $\mathbb{R}_+$, we define a relation $h_1(x) \asymp h_2(x)$ as $x \to \infty$ by $\limsup_{x \to \infty} h_2(x)/h_1(x) < \infty$ and $\liminf_{x \to \infty} h_2(x)/h_1(x) > 0$. For a positive increasing function $h(x)$ on $(0, \infty)$, we define the inverse function $h^{-1}(x)$ as

$$h^{-1}(x) := \sup\{y : h(y) < x\}$$

with the understanding that $\sup \varnothing = 0$. Let $\alpha := \log a$. We define two classes $\mathcal{G}$ and $\Phi$ of functions on $\mathbb{R}_+$, depending on $a$, as

(1.3) $\quad \mathcal{G} := \left\{ g(x) : g(x) > 0 \text{ on } \mathbb{R}_+ \text{ and } \limsup_{n \to \infty} g(n+1)/g(n) < a \right\}$

and

(1.4)
$$\begin{aligned}
\Phi := \{ \phi : \phi(t) = t^\alpha g(|\log t|) \text{ on } (0, \infty) \text{ with } g(x) \in \mathcal{G}, \\
\text{and } \phi(t) \text{ is positive} \\
\text{and increasing on } (0, \delta_1) \text{ with some } \delta_1 > 0 \\
\text{satisfying } \phi(0) := \phi(0+) = 0 \}.
\end{aligned}$$

Note that we assume the monotone property for $\phi \in \Phi$ but do not for $g \in \mathcal{G}$. For a nonnegative function $\phi$ on $(0, \infty)$, which is positive and increasing on $(0, \delta_1)$ for some $\delta_1 > 0$ with $\phi(0) := \phi(0+) = 0$, the $\phi$-Hausdorff measure $\phi$-$H(E)$ of a Borel set $E$ in the metric space $(\mathbf{I}, d)$ is defined by

(1.5) $\quad \phi\text{-}H(E) := \liminf_{\delta \to 0+} \left\{ \sum_{n=1}^\infty \phi(|D_n|) : E \subset \bigcup_{n=1}^\infty D_n, |D_n| \leq \delta \right\},$

where $|D_n|$ denotes the diameter of the set $D_n \in \mathcal{B}(\mathbf{I})$. We can take the set $D_n$ as a closed ball in the definition (1.5), since $(\mathbf{I}, d)$ is an ultrametric space. Under the single assumption that $a > 1$, both the Hausdorff and packing dimensions of $\partial \mathbf{T}$ are $\alpha$ almost surely on $\{\partial \mathbf{T} \neq \varnothing\}$. See, for Hausdorff dimension, [12, 17, 18] and [26]; for packing dimension, [4, 22] and [46]. See also [20] for an extension of Hawkes's result. A $\phi$-Hausdorff measure is called an *exact Hausdorff measure* for $\partial \mathbf{T}$ if $0 < \phi\text{-}H(\partial \mathbf{T}) < \infty$ a.s. on $\{\partial \mathbf{T} \neq \varnothing\}$. Denote by $\mu$ the *branching measure* on the boundary $\partial \mathbf{T}$. An exact Hausdorff measure is called absolutely continuous (with respect to the branching measure $\mu$) if $\phi\text{-}H(A) = 0$ a.s. provided that $\mu(A) = 0$ for a Borel set $A \subset \partial \mathbf{T}$. A precise definition of the branching measure $\mu$ is given in Section 2. The random sequence $\{Y(n), n \leq 0\}$ on an extended probability space $(\Omega \times \mathbf{I}, \mathcal{F} \times \mathcal{B}(\mathbf{I}), Q)$ will be defined by (2.6) and (2.7) in Section 2. There can be more than one exact Hausdorff measure. Our main results are as follows.

We define an integrated function $K(x)$ of a tail probability on $\mathbb{R}_+$ as

(1.6) $$K(x) := \int_x^\infty P(N > u) \, du.$$

We define an exceptional set $\Delta$ in $\partial \mathbf{T}$ and a condition $(G_\Delta)$ for $g \in \mathcal{G}$ as follows:

(1.7) $$\Delta := \left\{ \mathbf{i} \in \partial \mathbf{T} : \lim_{n \to \infty} \frac{W_{\mathbf{i}|n}}{g(n)} = 0 \right\},$$

where $W_{\mathbf{i}}$ is defined by (2.3) below and

$$(G_\Delta) \quad \limsup_{n\to\infty}\left(\sum_{k=0}^{n-1} Q(Y(0)-Y(-1) > \delta_0 g(k)) - \log g(n)\right) = \infty$$

for some $\delta_0 > 0$.

The relation between the set $\Delta$ and the condition $(G_\Delta)$ is found in Lemma 4.3 below.

THEOREM 1.1. *Suppose that $K(x) \notin \mathcal{D}$. Then there exists an absolutely continuous exact $\phi$-Hausdorff measure for $\partial \mathbf{T}$ with $\phi(t) = t^\alpha g(|\log t|) \in \Phi$. It satisfies that*

$$(1.8) \quad \phi\text{-}H(\partial \mathbf{T}) = C_\phi W \quad a.s. \ on \ \{\partial \mathbf{T} \neq \varnothing\},$$

*where the positive constant $C_\phi$ is determined by*

$$(1.9) \quad \sum_{n=0}^\infty Q(Y(0) - Y(-\ell) > \delta g(n)) \begin{cases} = \infty, & \exists \ell = \ell(\delta) \geq 1, \\ & \text{for } 0 < \delta < C_\phi^{-1}, \\ < \infty, & \forall \ell \geq 1, \ \text{for } \delta > C_\phi^{-1}. \end{cases}$$

*Moreover, it is represented, for $A \in \mathcal{B}(\mathbf{I})$ satisfying $A \subset \partial \mathbf{T}$, as*

$$(1.10) \quad \phi\text{-}H(A) = C_\phi \mu(A) \quad a.s.$$

REMARK 1.1. A concrete but not simple example of $\phi \in \Phi$ for an exact $\phi$-Hausdorff measure for $\partial \mathbf{T}$ is found in the proof of the above theorem in Section 4. The condition $K(x) \notin \mathcal{D}$ is equivalent to $P(W > x) \notin \mathcal{D}$, but not to $P(N > x) \notin \mathcal{D}$. See Lemma 2.5.

THEOREM 1.2. *Suppose that $K(x) \in \mathcal{D}$. Let $\phi \in \Phi$ with $\phi(t) = t^\alpha g(|\log t|)$ and let $\Delta'$ be an arbitrary Borel set in $\partial \mathbf{T}$ with $\mu(\Delta') = 0$ a.s. Then there is no absolutely continuous exact Hausdorff measure for $\partial \mathbf{T}$. More precisely, we have the following:*

(i) *If $\sum_{n=0}^\infty K(g(n)) < \infty$, then $\phi\text{-}H(\partial \mathbf{T} \setminus \Delta') = \phi\text{-}H(\partial \mathbf{T}) = \infty$ a.s. on $\{\partial \mathbf{T} \neq \varnothing\}$.*
(ii) *If $\sum_{n=0}^\infty K(g(n)) = \infty$, then $\phi\text{-}H(\partial \mathbf{T} \setminus \Delta) = 0$ with $\mu(\Delta) = 0$ a.s.*
(iii) *If $\sum_{n=0}^\infty K(g(n)) = \infty$ and*

$$\limsup_{\delta \to 0+} \limsup_{n\to\infty} \sum_{k=0}^n K(\delta g(k))/\log(g(n) \vee e) = \infty,$$

*then $\phi\text{-}H(\partial \mathbf{T}) = 0$ a.s.*

REMARK 1.2. In the case where $K(x) \in \mathcal{D}$, it is still hard to answer whether there exists an exact Hausdorff measure for $\partial \mathbf{T}$. However, our results say the following. Suppose that $K(x) \in \mathcal{D}$ and there exists an exact $\phi$-Hausdorff measure $\phi\text{-}H$ for $\partial \mathbf{T}$ with $\phi \in \Phi$. Then it is singular with respect to the branching measure $\mu$ and satisfies that

$$0 < \phi\text{-}H(\Delta) < \infty \quad \text{and} \quad \phi\text{-}H(\partial \mathbf{T} \setminus \Delta) = 0 \qquad \text{with } \mu(\Delta) = 0 \text{ a.s.}$$

Further it satisfies that $\lim_{n \to \infty} g(n) = \infty$ and that

$$\sum_{n=0}^{\infty} K(g(n)) = \infty \quad \text{and} \quad \limsup_{\delta \to 0+} \limsup_{n \to \infty} \sum_{k=0}^{n} K(\delta g(k)) / \log g(n) < \infty.$$

COROLLARY 1.1. (i) *There exists an absolutely continuous exact $\phi$-Hausdorff measure for $\partial \mathbf{T}$ with $\phi \in \Phi$ if and only if $K(x) \notin \mathcal{D}$, that is, $P(W > x) \notin \mathcal{D}$. It is also equivalent to that there exists $\phi \in \Phi$ such that $0 < \phi\text{-}H(\partial \mathbf{T} \setminus \Delta) < \infty$ with $\mu(\Delta) = 0$ a.s. on $\{\partial \mathbf{T} \neq \varnothing\}$.*

(ii) *An exact $\phi$-Hausdorff measure for $\partial \mathbf{T}$ with $\phi \in \Phi$ is absolutely continuous if and only if $K(x) \notin \mathcal{D}$ and $\phi\text{-}H(\Delta) = 0$ a.s. If an exact $\phi$-Hausdorff measure with $\phi \in \Phi$ satisfies the condition ($G_\Delta$) for its g function, then it is absolutely continuous.*

REMARK 1.3. There is symmetry between the problem on the exact Hausdorff measure on the boundary of the tree and that on the exact packing measure. The former is related to the right tail behavior of the distribution of $W$ and the latter is to the left tail behavior. The existence of an exact packing measure for the tree is determined by the nondominated variation of the left tail of the distribution of $W$, namely, $p_0 = p_1 = 0$. In addition, the exact packing measure is explicitly given and absolutely continuous with respect to the branching measure. See [46]. Moreover, $t^\alpha$-Hausdorff measure (so-called $\alpha$-Hausdorff measure) of $\partial \mathbf{T}$ is almost surely zero. On the other hand, $t^\alpha$-packing measure (so-called $\alpha$-packing measure) of $\partial \mathbf{T}$ is almost surely infinity on $\{\partial \mathbf{T} \neq \varnothing\}$. An analogous symmetry is already found in the results of [45].

We add the three theorems which were discussed by Liu [21]. Unfortunately, his proofs of Theorems 1.3 and 1.4 contain serious gaps as was already pointed out by Watanabe [46]. Theorems 1.3 and 1.4 correspond to Theorems 2 and 3 in [21]. See Remark 1.4. The results of Bingham [8], Watanabe [45] and Bingham and Doney [9] on the right tail behavior of the martingale limit $W$ are used in the proof of the following theorems, respectively.

THEOREM 1.3. *Let $M := \sup\{n \geq 0 : p_n > 0\}$. Suppose that $1 < M < \infty$. Define $\gamma$ and $\phi_0$ as $\gamma := \log M / \log a$ and $\phi_0(t) := t^\alpha (\log|\log t|)^{(\gamma-1)/\gamma}$. Then we have $\gamma > 1$ and*

$$(1.11) \qquad \phi_0\text{-}H(\partial \mathbf{T}) = \tau^{(\gamma-1)/\gamma} W \qquad \text{a.s. on } \{\partial \mathbf{T} \neq \varnothing\},$$

*where $\tau$ is a positive constant determined by*

$$(1.12) \qquad E \exp(\delta W^{\gamma/(\gamma-1)}) \begin{cases} < \infty, & \text{for } 0 < \delta < \tau, \\ = \infty, & \text{for } \delta > \tau. \end{cases}$$

THEOREM 1.4. *Let $S_0 := \sup\{s > 0 : f(s) < \infty\}$. Suppose that $1 < S_0 < \infty$. Define $\phi_1$ as $\phi_1(t) := t^\alpha \log|\log t|$. Then we have*

$$(1.13) \qquad \phi_1\text{-}H(\partial \mathbf{T}) = \sigma W \qquad \text{a.s. on } \{\partial \mathbf{T} \neq \varnothing\},$$

*where $\sigma$ is a positive constant given by*

$$(1.14) \qquad \sigma := \lim_{n \to \infty} a^{n+1}((f_n)^{-1}(S_0) - 1),$$

*where $(f_n)^{-1}(s)$ is the inverse function of $f_n(s)$.*

THEOREM 1.5. *Let $\theta_0 := \sup\{\theta \geq 1 : E(N^\theta) < \infty\}$. Suppose that $1 < \theta_0 < \infty$. Define $b_0$ and $\psi_b$ as $b_0 := 1/(\theta_0 - 1)$ and $\psi_b(t) := t^\alpha |\log t|^b$ for $-\infty < b < \infty$. Then we have, a.s. on $\{\partial \mathbf{T} \neq \varnothing\}$,*

$$(1.15) \qquad \psi_b\text{-}H(\partial \mathbf{T}) = \begin{cases} 0, & \text{for } b < b_0, \\ \infty, & \text{for } b > b_0. \end{cases}$$

*Moreover, if $E(N^{\theta_0}) < \infty$, then $\psi_{b_0}\text{-}H(\partial \mathbf{T}) = \infty$ a.s. on $\{\partial \mathbf{T} \neq \varnothing\}$. On the other hand, if $E(N^{\theta_0}) = \infty$, then $\psi_{b_0}\text{-}H(\partial \mathbf{T} \setminus \Delta) = 0$ with $\mu(\Delta) = 0$ a.s.*

Finally we present a resolution for a conjecture of Hawkes [17]. However, it should be noted that no necessary and sufficient condition for the relation (1.16) below in terms of the offspring distribution of the branching process is known up to now.

THEOREM 1.6. *Let $R(x) := x^b \ell(x)$ be a positive and increasing function on $(0, \infty)$ with $0 < b \leq 1$ and slowly varying $\ell(x)$ as $x \to \infty$. Define $\phi_2$ as $\phi_2(t) := t^\alpha R^{-1}(\log(e \vee |\log t|))$. Suppose that $\phi_2 \in \Phi$ and*

$$(1.16) \qquad -\log P(W > x) \asymp R(x) \qquad \text{as } x \to \infty.$$

*Then we have*

$$(1.17) \qquad \phi_2\text{-}H(\partial \mathbf{T}) = \xi_R W \qquad \text{a.s. on } \{\partial \mathbf{T} \neq \varnothing\},$$

*where $\xi_R$ is a positive constant determined by*

$$(1.18) \qquad E \exp(R(\delta W)) \begin{cases} < \infty, & \text{for } 0 < \delta < \xi_R, \\ = \infty, & \text{for } \delta > \xi_R. \end{cases}$$

REMARK 1.4. Each one-half of the proofs of Theorems 1.3 and 1.4 by Liu [21] has a serious gap. Namely, the first equality on line 7 on page 535 and the reason for the inequality on line 11 on page 536 of [21] are not true, respectively. The first assertion of Theorem 1.5 was conjectured by Liu [21]. All the exact Hausdorff measures $\phi_j$-$H$ ($j = 0, 1, 2$) in Theorems 1.3, 1.4 and 1.6 are absolutely continuous because they satisfy the condition ($G_\Delta$) for their $g$ functions. Thus they satisfy (1.10) with $C_{\phi_0} = \tau^{(\gamma-1)/\gamma}$, $C_{\phi_1} = \sigma$ and $C_{\phi_2} = \xi_R$, respectively satisfying $\phi_j$-$H(\Delta) = 0$ a.s. ($j = 0, 1, 2$). The constants $\tau$ and $\sigma$ are explained more precisely in Lemmas 2.4 and 2.5 and Remarks 2.1 and 2.2. The existence and positivity of $\xi_R$ are trivial, but the method of the numerical calculation of its value is not known. We do not know whether $\psi_{b_0}$-$H(\partial \mathbf{T}) = 0$ a.s., that is, whether $\psi_{b_0}$-$H(\Delta) = 0$ a.s., in the case where $E(N^{\theta_0}) = \infty$.

The organization of this paper is as follows. In Section 2, we review some useful results on the distribution of $W$ and those on shift self-similar additive random sequences, and give some preliminary results. In Section 3, we give limit theorems of "limsup" type for a shift self-similar additive random sequence $\{Y(n)\}$. In Section 4, we prove the main theorems and the four additional theorems.

**2. Preliminaries.** Let $c > 1$. An $\mathbb{R}^d$-valued random sequence $\{X(n), n \in \mathbb{Z}\}$ on a probability space $(\widetilde{\Omega}, \widetilde{\mathcal{F}}, \widetilde{P})$ is called a *shift c-self-similar additive random sequence* if the following two conditions are satisfied:

(1) The sequence $\{X(n), n \in \mathbb{Z}\}$ has shift $c$-self-similarity, that is,

$$\{X(n+1), n \in \mathbb{Z}\} \stackrel{\mathrm{d}}{=} \{cX(n), n \in \mathbb{Z}\},$$

where the symbol $=^{\mathrm{d}}$ stands for equality in the finite-dimensional distributions.

(2) The sequence $\{X(n), n \in \mathbb{Z}\}$ has independent increments (or additivity), that is, for every $n \in \mathbb{Z}$, $\{X(k), k \leq n\}$ and $X(n+1) - X(n)$ are independent.

The definition for an $\mathbb{R}^d$-valued random sequence $\{X(n), n \leq 0\}$ to be a shift $c$-self-similar additive random sequence is similar. That is, the sequence $\{X(n), n \leq 0\}$ is called a shift $c$-self-similar additive random sequence if $\{X(n+1), n \leq -1\} =^{\mathrm{d}} \{cX(n), n \leq -1\}$ and, for every $n \leq -1$, $\{X(k), k \leq n\}$ and $X(n+1) - X(n)$ are independent. Note that shift self-similarity does not imply the usual self-similarity. We denote by $\widehat{\eta}$ the characteristic function of a probability distribution $\eta$ on $\mathbb{R}^d$. Let $0 < b < 1$. A probability distribution $\eta$ on $\mathbb{R}^d$ is said to be *b-decomposable* if there exists a probability distribution $\rho$ on $\mathbb{R}^d$ such that

$$\widehat{\eta}(z) = \widehat{\eta}(bz)\widehat{\rho}(z). \tag{2.1}$$

For example, semistable distributions and homogeneous self-similar measures such as Bernoulli convolutions are $b$-decomposable for some $b \in (0,1)$. In the case where $\rho$ is infinitely divisible, $\eta$ is also infinitely divisible and is called *semi-self-decomposable*. The equality (2.1) is equivalent to

$$\widehat{\eta}(z) = \prod_{n=0}^{\infty} \widehat{\rho}(b^n z) \tag{2.2}$$

which is convergent if and only if $\int_{\mathbb{R}^d} \log(1+|x|)\rho(dx) < \infty$. The distribution $\rho$ in (2.1) is not necessarily uniquely determined by the distribution $\eta$. It is uniquely determined by the distribution $\eta$ in case the support of $\eta$ is contained in $\mathbb{R}_+^d$. See [11, 25, 27, 29] and [30]. Absolute continuity of $b$-decomposable distributions is very difficult and is related to *Peres–Solomyak numbers* (PS numbers, for short) and *Pisot–Vijayaraghavan numbers* (PV numbers, for short). See [42, 43] and [47]. Let $\{V_j\}$ be i.i.d. random variables with the distribution $\rho$ in (2.1). Then the sequence $\{X(n), n \in \mathbb{Z}\}$ defined by $X(n) = \sum_{j=-\infty}^{n} c^j V_j$ with $c = b^{-1}$ is a shift $c$-self-similar additive random sequence with $\eta$ being the distribution of $X(0)$. For instance, let $\{T_n\}$ be a first exit time from an equilateral triangle with side $2^n$ of a Brownian motion $\{B(t)\}$ on the extended Sierpinski gasket, starting at the origin. Then $\{T_n, n \in \mathbb{Z}\}$ is a shift 5-self-similar additive random sequence and, for each $n \in \mathbb{Z}$ the distribution of $T_n$ is semi-self-decomposable and absolutely continuous with infinitely differentiable density. Certain limit theorems of "limsup" type and "liminf" type for $\{T_n\}$ are equivalent to the laws of the iterated logarithm of "liminf" type and "limsup" type for $\{B(t)\}$, respectively. In particular, the constants in the laws of the iterated logarithm for $\{B(t)\}$ are unknown up to now, but their upper and lower bounds are explicitly obtained together with their candidates by using analogous constants in the limit laws for $\{T_n\}$. The same kinds of results are also true for Brownian motions on many nested fractals other than Sierpinski gasket. See [42, 45] and Remark 2.2.

The author [44], motivated by the results of Sato [36], introduced and characterized shift self-similar additive random sequences, and studied in detail their transience and rate of growth. Further, he found in [45, 46] two important examples of them in relation to general supercritical Galton–Watson branching processes. See Theorem 1.1 of [45] and Theorem 2.1 of [46]. They are closely related to self-similar or semi-self-similar additive processes. See [28, 31, 37] and [41]. In general, finite-dimensional distributions of the shift self-similar additive random sequence $\{X(n)\}$ are determined by the distribution of $X(0) - X(-1)$ but not always by that of $X(0)$. We shall use the following increasing case.

LEMMA 2.1 (Theorem 2.1 of [44]). *Let $c > 1$.*

(i) *Suppose that $\{X(n), n \leq 0\}$ is an increasing shift c-self-similar additive random sequence. Then the distribution of $X(n)$ is $c^{-1}$-decomposable on $\mathbb{R}_+$ for $n \leq 0$ and $\lim_{n \to \infty} X(-n) = 0$ a.s. There is a unique in law increasing shift c-self-similar additive random sequence $\{\widetilde{X}(n), n \in \mathbb{Z}\}$ such that*

$$\{X(n), n \leq 0\} \stackrel{\mathrm{d}}{=} \{\widetilde{X}(n), n \leq 0\}.$$

(ii) *Conversely, for any $c^{-1}$-decomposable distribution $\eta$ on $\mathbb{R}_+$, there exists a unique in law increasing shift c-self-similar additive random sequence $\{X(n), n \leq 0\}$ with the distribution of $X(0)$ being $\eta$.*

Let $\mathbf{T}$ be a Galton–Watson tree on $(\Omega, \mathcal{F}, P)$ with $f(s)$ being the offspring p.g.f. We continue to use the notation and the assumptions in Section 1. In particular, the random variable $W$ is defined by (1.2). We define a shifted tree $\mathbf{T_i}$ of $\mathbf{T}$ for $\mathbf{i} \in \mathbf{U}$ by the following two conditions:

(1) $\varnothing \in \mathbf{T_i}$.
(2) Let $\mathbf{j} \in \mathbf{T_i}$ and $i \in \mathbb{Z}_+$. Then $\mathbf{j} * i \in \mathbf{T_i}$ if and only if $0 \leq i \leq N_{\mathbf{i}*\mathbf{j}} - 1$.

We define an $\mathbb{R}_+$-valued random variable $W_\mathbf{i}$ for $\mathbf{i} \in \mathbf{U}$ as

$$(2.3) \qquad W_\mathbf{i} := \lim_{n \to \infty} \frac{\mathrm{Card}\{\mathbf{j} \in \mathbf{T_i} : |\mathbf{j}| = n\}}{a^n}.$$

The limit $W_\mathbf{i}$ exists almost surely. It satisfies that

$$(2.4) \qquad W_\mathbf{i} = \frac{1}{a} \sum_{j=0}^{N_\mathbf{i} - 1} W_{\mathbf{i}*j} \qquad \text{for } \mathbf{i} \in \mathbf{U},$$

with the understanding that $\sum_{j=0}^{-1} = 0$. The distribution of $W_\mathbf{i}$ is the same as that of $W$ for $\mathbf{i} \in \mathbf{U}$. Moreover, $W_\mathbf{i}$ and $W_\mathbf{j}$ are independent, whenever neither $\mathbf{i} \leq \mathbf{j}$ nor $\mathbf{j} \leq \mathbf{i}$. We define a closed ball $B_\mathbf{i}$ in $\mathbf{I}$ and its diameter $|B_\mathbf{i}|$ for $\mathbf{i} \in \mathbf{U}$ as

$$B_\mathbf{i} := \{\mathbf{j} \in \mathbf{I} : \mathbf{i} \leq \mathbf{j}\} \quad \text{and} \quad |B_\mathbf{i}| = e^{-|\mathbf{i}|}.$$

Note that $\{W_\mathbf{i} = 0\} = \{\partial \mathbf{T} \cap B_\mathbf{i} = \varnothing\}$ for $\mathbf{i} \in \mathbf{T}$ up to a probability zero set. A subset $\Gamma$ of $\mathbf{T}$ is called a *cutset* for a subset $A$ of the boundary $\partial \mathbf{T}$ if the following three conditions are satisfied:

(1) Neither $\mathbf{i} \leq \mathbf{j}$ nor $\mathbf{j} \leq \mathbf{i}$ whenever $\mathbf{i}, \mathbf{j} \in \Gamma$ and $\mathbf{i} \neq \mathbf{j}$.
(2) $A \subset \bigcup_{\mathbf{i} \in \Gamma} B_\mathbf{i}$.
(3) For $\mathbf{i} \in \Gamma$, $B_\mathbf{i} \cap A \neq \varnothing$.

Let $\Gamma$ be any cutset for $\partial \mathbf{T} \cap B_\mathbf{i}$. We see from (2.4) that

$$(2.5) \qquad |B_\mathbf{i}|^\alpha W_\mathbf{i} = \sum_{\mathbf{j} \in \Gamma} |B_\mathbf{j}|^\alpha W_\mathbf{j}.$$

We define a finite measure $\mu = \mu_\omega$ on $(\mathbf{I}, \mathcal{B}(\mathbf{I}))$ by

$$\mu(B_\mathbf{i}) := \begin{cases} a^{-|\mathbf{i}|} W_\mathbf{i}, & \text{for } \mathbf{i} \in \mathbf{T}, \\ 0, & \text{for } \mathbf{i} \in \mathbf{U} \setminus \mathbf{T}. \end{cases}$$

Note that $\mu$ is determined uniquely on $(\mathbf{I}, \mathcal{B}(\mathbf{I}))$ for each $\omega \in \Omega$ and the support of $\mu$ is contained in $\partial \mathbf{T}$ almost surely. The measure $\mu$ is called the branching measure for the tree $\mathbf{T}$. See [21] as to the above assertions. We define a probability space $(\Omega \times \mathbf{I}, \mathcal{F} \times \mathcal{B}(\mathbf{I}), Q)$ by assigning $Q(A)$, for $A \in \mathcal{F} \times \mathcal{B}(\mathbf{I})$,

$$(2.6) \qquad Q(A) := E\left(\int \mathbb{1}_A(\omega, \mathbf{i}) \mu_\omega(d\mathbf{i})\right),$$

where $\mathbb{1}_A$ stands for the indicator function of a set $A$. We denote by $E_Q$ the expectation under the probability measure $Q$. We define a random sequence $\{Y(n), n \leq 0\}$ by

$$(2.7) \qquad Y(-n) := \mu(B_{\mathbf{i}|n}) \qquad \text{for } n \in \mathbb{Z}_+ \text{ and } \mathbf{i} \in \mathbf{I}.$$

The shift self-similarity of $\{Y(n)\}$ was suggested by Hawkes [17], but not so was the additivity in the following lemma.

LEMMA 2.2 (Theorem 2.1 of [46]). *The sequence $\{Y(n), n \leq 0\}$ is an increasing $\mathbb{R}_+$-valued shift $a$-self-similar additive random sequence on the probability space $(\Omega \times \mathbf{I}, \mathcal{F} \times \mathcal{B}(\mathbf{I}), Q)$. In particular, we have $Q(Y(0) \leq x) = E(W \mathbb{1}_{\{W \leq x\}})$ and the distribution of $Y(0)$ under $Q$ is $a^{-1}$-decomposable.*

For a distribution $\eta$ on $\mathbb{R}_+$, we denote the tail by $\bar{\eta}(x)$, that is, $\bar{\eta}(x) = \eta(x, \infty)$ for $x \geq 0$. We denote the convolution of distributions $\eta$ and $\nu$ by $\eta * \nu$ and denote the $n$th convolution power of $\eta$ by $\eta^{n*}$ with the understanding that $\eta^{0*}(dx) = \delta_0(dx)$, namely, the probability distribution concentrated at 0. Denote the distribution of $W$ under the probability measure $P$ by $\nu_W$. The characteristic function of the distribution $\nu_W$ satisfies Poincaré's equation, that is,

$$(2.8) \qquad \widehat{\nu}_W(z) = f(\widehat{\nu}_W(z/a)).$$

Define the p.g.f. $\widetilde{f}(s) = \sum_{n=0}^\infty \widetilde{p}_n s^n$ and the distribution $\eta_W$ by

$$(2.9) \quad \widetilde{f}(s) := \frac{f(q + (1-q)s) - q}{1 - q} \quad \text{and} \quad \eta_W(dx) := \frac{\nu_W(dx) - q\delta_0(dx)}{1 - q}.$$

Note that $\widetilde{p}_0 = 0$ and $\bar{\eta}_W(x) = (1-q)^{-1} \bar{\nu}_W(x) > 0$ for any $x > 0$ and that $\nu_W = \eta_W$ if and only if $p_0 = 0$. Define the distributions $\eta'_W$ and $\rho_W$ by

$$(2.10) \quad \eta'_W(dx) := \eta_W(a\,dx) \quad \text{and} \quad \rho_W(dx) := \sum_{n=1}^\infty \widetilde{p}_n (\eta'_W)^{(n-1)*}(dx).$$

Then we obtain from (2.8) and (2.9) that

(2.11)    $\widehat{\eta}_W(z) = \widetilde{f}(\widehat{\eta}_W(z/a))$   and   $\eta_W(dx) = \eta'_W * \rho_W(dx).$

That is, $\eta_W$ is $a^{-1}$-decomposable. We denote the distributions of $Y(0)$ and $Y(0) - Y(-1)$ under the probability measure $Q$ by $\eta_Y$ and $\rho_Y$. Then we have by (2.2) and Lemma 2.2

(2.12)    $\eta_Y(dx) = \eta'_Y * \rho_Y(dx)$   and   $\int_1^\infty (\log x)\rho_Y(dx) < \infty,$

where $\eta'_Y(dx) := \eta_Y(a\,dx)$.

Bingham [8] proved the following lemma by using a Tauberian theorem of exponential type of Kasahara [19] and Theorem 3.4 of [16].

LEMMA 2.3 ((11) of [8]).   *Let $M := \sup\{n \geq 0 : p_n > 0\}$. Suppose that $1 < M < \infty$. Let $\gamma = \log M/\log a$. Then we have $\gamma > 1$ and*

(2.13)    $-\log(\bar{\nu}_W(x)) \asymp x^{\gamma/(\gamma-1)}$   *as $x \to \infty$.*

*Moreover, we see that there exists a positive constant $\tau$ determined by (1.12) and (1.12) is clearly equivalent to*

(2.14)    $E(W \exp(\delta W^{\gamma/(\gamma-1)})) \begin{cases} < \infty, & \text{for } 0 < \delta < \tau, \\ = \infty, & \text{for } \delta > \tau. \end{cases}$

REMARK 2.1.   The positive constant $\tau$ is represented by Liu [22] as

$$\tau = \liminf_{x\to\infty}\{-x^{\gamma/(1-\gamma)} \log P(W > x)\}.$$

Numerical calculation of the value of $\tau$ is very difficult. See [7] and also [5, 6].

LEMMA 2.4 (Theorem 2.2 of [45]).   *Let $S_0 := \sup\{s > 0 : f(s) < \infty\}$. Suppose that $1 < S_0 < \infty$. Then there exists a positive constant $\sigma$ such that*

(2.15)    $E \exp(tW) \begin{cases} < \infty, & \text{for } 0 < t < \sigma, \\ = \infty, & \text{for } t > \sigma. \end{cases}$

*Further, it is represented as (1.14) and (2.15) is obviously equivalent to*

(2.16)    $E(W \exp(tW)) \begin{cases} < \infty, & \text{for } 0 < t < \sigma, \\ = \infty, & \text{for } t > \sigma. \end{cases}$

REMARK 2.2.   Suppose that, for some positive integer $k$,

$$f(s) = \frac{s}{(a - (a-1)s^k)^{1/k}}.$$

Then we see that

$$f_n(s) = \frac{s}{(a^n - (a^n - 1)s^k)^{1/k}}$$

and $\widehat{\nu}_W(z) = (1 - ikz)^{-1/k}$, that is, $\nu_W$ is the gamma distribution with parameter $1/k$ and thereby $\sigma = 1/k$. See [16] or [8]. The existence and the positivity of $\sigma$ are also found in [22] but its representation (1.14) is not therein. It should be noted that the constant $\sigma$ can be numerically calculated by using (1.14). In some cases, the constant $\sigma$ has a natural relation with the constant of the law of the iterated logarithm of Brownian motions on some fractals. See [2, 13] and [45]. The Brownian motion on the Sierpinski gasket is related to the case where $f(s) = s^2/(4-3s)$ with $a = 5$ and $S_0 = 4/3$, and $\sigma$ is computed numerically by using (1.14) as $\sigma = 1.318\cdots$.

The author [45] proved the following lemma by using a Tauberian theorem of [41]. It is the most difficult and critical fact in this paper. The regularly varying case was already known by Bingham and Doney [9] and de Meyer [32] in a stronger sense. Recall the definition (1.6) of $K(x)$.

LEMMA 2.5 (Theorem 2.3 of [45]).  (i) $\bar{\nu}_W(x) \in \mathcal{D}$ if and only if $K(x) \in \mathcal{D}$.

(ii) If $K(x) \in \mathcal{D}$, then

(2.17) $$x\bar{\nu}_W(x) \asymp K(x) \qquad as\ x \to \infty.$$

REMARK 2.3.  We can prove as in the proof of Lemma 4.1 of [45] by using (2.11) that, for some $c_1 > 0$, $c_1 \bar{\eta}_W(ax) \leq \bar{\rho}_W(x) \leq \bar{\eta}_W(x)$ for $x \geq 0$. It follows that $\bar{\rho}_W(x) \in \mathcal{D}$ if and only if $\bar{\eta}_W(x) \in \mathcal{D}$ and that if $\bar{\eta}_W(x) \in \mathcal{D}$, then $\bar{\rho}_W(x) \asymp \bar{\eta}_W(x)$. Thus it follows from the above lemma that if $K(x) \in \mathcal{D}$, then $x\bar{\rho}_W(x) \asymp x\bar{\eta}_W(x) \asymp x\bar{\nu}_W(x) \asymp K(x)$ as $x \to \infty$.

LEMMA 2.6.  Suppose that $\eta$ is $a^{-1}$-decomposable on $\mathbb{R}_+$ such that $\eta(dx) = \eta' * \rho(dx)$ with $\eta'(dx) := \eta(a\, dx)$ and $\int_1^\infty (\log x)\rho(dx) < \infty$.

(i) We have, for some $\varepsilon \in (0,1)$

(2.18) $\quad \bar{\rho}(x)(1 - \bar{\eta}(ax)) \leq \bar{\eta}(x) - \bar{\eta}(ax) \leq 2\bar{\rho}(\varepsilon x) \qquad for\ x > 0.$

(ii) $\bar{\rho}(x) \in \mathcal{D}$ if and only if $\bar{\eta}(x) - \bar{\eta}(ax) \in \mathcal{D}$.
(iii) If $\bar{\rho}(x) \in \mathcal{D}$, then

(2.19) $$\bar{\rho}(x) \asymp \bar{\eta}(x) - \bar{\eta}(ax) \qquad as\ x \to \infty.$$

PROOF.  Let $x > 0$. It follows from (2.1) that

$$\bar{\eta}(x) = \bar{\rho}(x) + \int_{0-}^{x+} \bar{\eta}(a(x-y))\rho(dy)$$
$$\geq \bar{\rho}(x) + \bar{\eta}(ax)(1 - \bar{\rho}(x)),$$

that is,

(2.20) $$\bar{\rho}(x)(1 - \bar{\eta}(ax)) \leq \bar{\eta}(x) - \bar{\eta}(ax).$$

On the other hand, we obtain from (2.1) that, for $0 < \varepsilon < 1$,

$$\bar{\eta}(x) = \bar{\eta}(ax) + \int_{0-}^{x+} \bar{\rho}(x-y)\eta(a\,dy)$$

$$\leq \bar{\eta}(ax) + \int_{(1-\varepsilon)x+}^{x+} \eta(a\,dy) + \int_{0-}^{(1-\varepsilon)x+} \bar{\rho}(\varepsilon x)\eta(a\,dy)$$

$$\leq \bar{\eta}((1-\varepsilon)ax) + \bar{\rho}(\varepsilon x),$$

namely,

(2.21) $$\bar{\eta}(x) - \bar{\eta}((1-\varepsilon)ax) \leq \bar{\rho}(\varepsilon x).$$

Letting $\varepsilon$ satisfy $(1-\varepsilon)^2 a = 1$ and adding the following to (2.21):

$$\bar{\eta}((1-\varepsilon)ax) - \bar{\eta}((1-\varepsilon)^2 a^2 x) \leq \bar{\rho}(\varepsilon(1-\varepsilon)ax),$$

we see that

(2.22) $$\bar{\eta}(x) - \bar{\eta}(ax) \leq \bar{\rho}(\varepsilon x) + \bar{\rho}(\varepsilon(1-\varepsilon)ax) \leq 2\bar{\rho}(\varepsilon x).$$

Thus we have established (2.18) from (2.20) and (2.22). Assertions (ii) and (iii) are obvious from assertion (i). □

LEMMA 2.7. *(i)* $\bar{\rho}_Y(x) \in \mathcal{D}$ *if and only if* $K(x) \in \mathcal{D}$.
*(ii) If* $K(x) \in \mathcal{D}$, *then*

(2.23) $$\bar{\rho}_Y(x) \asymp K(x) \qquad \text{as } x \to \infty.$$

PROOF. It follows from Lemma 2.6 that $\bar{\rho}_Y(x) \in \mathcal{D}$ if and only if $\bar{\eta}_Y(x) - \bar{\eta}_Y(ax) \in \mathcal{D}$ and that if $\bar{\rho}_Y(x) \in \mathcal{D}$, then $\bar{\rho}_Y(x) \asymp \bar{\eta}_Y(x) - \bar{\eta}_Y(ax)$. In the same way, $\bar{\rho}_W(x) \in \mathcal{D}$ if and only if $\bar{\eta}_W(x) - \bar{\eta}_W(ax) \in \mathcal{D}$. Moreover, if $\bar{\rho}_W(x) \in \mathcal{D}$, then $\bar{\rho}_W(x) \asymp \bar{\eta}_W(x) - \bar{\eta}_W(ax)$. We see from Lemma 2.2 that

$$\bar{\eta}_Y(x) - \bar{\eta}_Y(ax) = \int_{x+}^{ax+} y\nu_W(dy)$$

(2.24) $$\asymp x(\bar{\nu}_W(x) - \bar{\nu}_W(ax))$$

$$\asymp x(\bar{\eta}_W(x) - \bar{\eta}_W(ax)).$$

Thus we conclude from Lemma 2.6 and Remark 2.3 that each of the six conditions $\bar{\rho}_Y(x) \in \mathcal{D}$, $\bar{\eta}_Y(x) - \bar{\eta}_Y(ax) \in \mathcal{D}$, $\bar{\eta}_W(x) - \bar{\eta}_W(ax) \in \mathcal{D}$, $\bar{\rho}_W(x) \in \mathcal{D}$, $\bar{\eta}_W(x) \in \mathcal{D}$ and $\bar{\nu}_W(x) \in \mathcal{D}$ is equivalent to $K(x) \in \mathcal{D}$. Further, if $K(x) \in \mathcal{D}$, then

$$\bar{\rho}_Y(x) \asymp \bar{\eta}_Y(x) - \bar{\eta}_Y(ax) \asymp x(\bar{\eta}_W(x) - \bar{\eta}_W(ax))$$

$$\asymp x\bar{\rho}_W(x) \asymp x\bar{\eta}_W(x) \asymp x\bar{\nu}_W(x) \asymp K(x).$$

Thus we have proved the lemma. □

LEMMA 2.8 (Theorem 5 of [9]). *Let $\theta > 1$. Then $E(W^\theta) < \infty$ if and only if $E(N^\theta) < \infty$.*

LEMMA 2.9. *Let $\theta > 0$. Then $E_Q((Y(0) - Y(-1))^\theta) < \infty$ if and only if $E(N^{\theta+1}) < \infty$.*

PROOF. Let $\theta > 0$. We obtain from (2.16) and (2.18) that
$$E_Q((Y(0) - Y(-1))^\theta) < \infty$$
$$\Leftrightarrow \int_{0-}^\infty x^\theta \rho_Y(dx) < \infty$$
$$\Leftrightarrow \int_0^\infty x^{\theta-1} \bar{\rho}_Y(x)\, dx < \infty$$
$$\Leftrightarrow \int_0^\infty x^{\theta-1} (\bar{\eta}_Y(x) - \bar{\eta}_Y(ax))\, dx < \infty.$$

By the same way we see from Lemma 2.8, Remark 2.3 and (2.18) that
$$E(N^{\theta+1}) < \infty \Leftrightarrow \int_{0-}^\infty x^{\theta+1} \nu_W(dx) < \infty$$
$$\Leftrightarrow \int_0^\infty x^\theta \bar{\nu}_W(x)\, dx < \infty$$
$$\Leftrightarrow \int_0^\infty x^\theta \bar{\eta}_W(x)\, dx < \infty$$
$$\Leftrightarrow \int_0^\infty x^\theta \bar{\rho}_W(x)\, dx < \infty$$
$$\Leftrightarrow \int_0^\infty x^\theta (\bar{\eta}_W(x) - \bar{\eta}_W(ax))\, dx < \infty.$$

Thus we conclude by (2.24) that $E_Q((Y(0) - Y(-1))^\theta) < \infty$ if and only if $E(N^{\theta+1}) < \infty$. □

LEMMA 2.10. *Let $\ell \in \mathbb{N}$. Then we have*
$$(2.25) \qquad \int_1^\infty Q(Y(0) - Y(-\ell) > x)\frac{dx}{x} < \infty.$$

PROOF. We consider on the probability space $(\Omega \times \mathbf{I}, \mathcal{F} \times \mathcal{B}(\mathbf{I}), Q)$. Since $Y(0) - Y(-\ell)$ and $Y(-\ell)$ are independent and the distributions of $Y(-\ell)$ and $a^{-\ell}Y(0)$ are the same, the distribution of $Y(0)$ is $a^{-\ell}$-decomposable. Hence the log-moment of $Y(0) - Y(-\ell)$ is finite by (2.2). That is,
$$\int_1^\infty Q(Y(0) - Y(-\ell) > x)\frac{dx}{x} = E_Q(\log((Y(0) - Y(-\ell)) \vee 1)) < \infty. \quad \square$$

**3. Limit theorems for $\{Y(n)\}$.** The author [46] studied the "liminf" type limit theorems for the sequence $\{Y(n)\}$. In this section, we discuss the "limsup" type limit theorems for the sequence $\{Y(n)\}$ by improving the results of [45]. Namely, we study the exact local dimension at typical $\mathbf{i} \in \partial \mathbf{T}$ of the branching measure $\mu$. Let $h(x)$ be positive measurable function on $(A, \infty)$ with $A \geq 0$. We say $h(x)$ is *submultiplicative* on $(A, \infty)$ if there is $c > 0$ such that $h(x+y) \leq ch(x)h(y)$ for all $x, y > A$. Further we say $h(x)$ is *quasi-submultiplicative* on $(A, \infty)$ if, for each $\varepsilon > 0$, there are $c_1, c_2 > 0$ such that $h(x+y) \leq c_1 h((1+\varepsilon)x) h(c_2 y)$ for all $x, y > A$. Obviously every submultiplicative function is quasi-submultiplicative but the converse is not true. For example, $(1 \vee x)^c$ with $c > 0$ is submultiplicative and $\exp(b_1 x^{b_2})$ with $b_1 > 0$ is submultiplicative for $0 < b_2 \leq 1$ and not so but quasi-submultiplicative for $1 < b_2 < \infty$ on $(0, \infty)$, respectively.

THEOREM 3.1. *Let $g \in \mathcal{G}$ and $C \in [0, \infty]$. We have*

$$\limsup_{n \to \infty} \frac{a^n Y(-n)}{g(n)} = C \qquad Q\text{-a.s.} \tag{3.1}$$

*if and only if*

$$\sum_{n=0}^{\infty} Q(Y(0) - Y(-\ell) > \delta g(n)) \begin{cases} = \infty, & \exists \ell = \ell(\delta) \geq 1, \ for \ 0 < \delta < C, \\ < \infty, & \forall \ell \geq 1, \ for \ \delta > C. \end{cases} \tag{3.2}$$

*Thus there is $C \in [0, \infty]$ satisfying (3.1) for each $g \in \mathcal{G}$.*

The proof of the above theorem is obvious from the following lemma.

LEMMA 3.1. *Let $g \in \mathcal{G}$ and $\ell \in \mathbb{N}$.*

(i) *If*

$$\sum_{n=0}^{\infty} Q(Y(0) - Y(-\ell) > g(n)) = \infty, \tag{3.3}$$

*then*

$$\limsup_{n \to \infty} \frac{a^n Y(-n)}{g(n)} \geq 1 \qquad Q\text{-a.s.} \tag{3.4}$$

(ii) *If*

$$\sum_{n=0}^{\infty} Q(Y(0) - Y(-\ell) > g(n)) < \infty, \tag{3.5}$$

*then*

$$\limsup_{n \to \infty} \frac{a^n Y(-n)}{g(n)} \leq \frac{1}{1 - (a_0/a)^\ell} \qquad Q\text{-a.s.,} \tag{3.6}$$

*where $a_0 := \limsup_{n \to \infty} g(n+1)/g(n) < a$ for $g \in \mathcal{G}$.*

PROOF. Let $g \in \mathcal{G}$. First we prove assertion (i). Suppose that (3.3) holds for some $\ell \geq 1$. Then it follows from the shift self-similarity that

$$\sum_{n=0}^{\infty} Q(Y(-n) - Y(-n-\ell) > a^{-n}g(n)) = \infty.$$

Thus there is $j_0$ with $0 \leq j_0 \leq \ell - 1$ such that

$$\sum_{n=0}^{\infty} Q(Y(-n\ell + j_0) - Y(-(n+1)\ell + j_0) > a^{-n\ell+j_0}g(n\ell - j_0)) = \infty.$$

Thanks to Borel–Cantelli's lemma, we have $Q$-a.s.

$$Y(-n\ell + j_0) - Y(-(n+1)\ell + j_0) > a^{-n\ell+j_0}g(n\ell - j_0) \quad \text{i.o.}$$

Here the abbreviation "i.o." stands for "infinitely often." Hence $Q$-a.s.

$$\frac{a^{n\ell-j_0}Y(-n\ell + j_0)}{g(n\ell - j_0)} > 1 \quad \text{i.o.}$$

Therefore we obtain (3.4). Next we prove assertion (ii). Suppose that (3.5) holds for some $\ell \geq 1$. Then there is a positive integer $n_0 = n_0(\omega)$ such that, for any $n \geq n_0$,

$$Y(-n) - Y(-n-\ell) \leq a^{-n}g(n) \quad Q\text{-a.s.}$$

Note from Lemma 2.1 that

$$\lim_{n \to \infty} Y(-n) = 0 \quad Q\text{-a.s.}$$

Hence we see that, for any $n \geq n_0$,

$$Y(-n) = \sum_{j=0}^{\infty} (Y(-n-j\ell) - Y(-n-(j+1)\ell))$$

$$\leq \sum_{j=0}^{\infty} a^{-n-j\ell}g(n+j\ell)$$

$$\leq a^{-n}g(n) \frac{1}{1 - ((a_0 + \varepsilon)/a)^{\ell}} \quad Q\text{-a.s.}$$

where the positive number $\varepsilon$ can be arbitrarily small when we take $n_0$ sufficiently large. Thus we have (3.6). □

THEOREM 3.2. Let $g \in \mathcal{G}$. Suppose that $K(x) \in \mathcal{D}$. If

$$\sum_{n=0}^{\infty} K(g(n)) = \infty \quad (resp. \ < \infty),$$

then

$$\limsup_{n \to \infty} \frac{a^n Y(-n)}{g(n)} = \infty \quad (resp. \ 0) \ Q\text{-a.s.}$$

PROOF. Suppose that $K(x) \in \mathcal{D}$. Then we see from Lemma 2.5 that

$$\sum_{n=0}^{\infty} Q(Y(0) - Y(-1) > \delta g(n)) = \infty \qquad (\text{resp.} < \infty) \text{ for any } \delta > 0,$$

if and only if $\sum_{n=0}^{\infty} K(g(n)) = \infty$ (resp. $< \infty$). Thus we obtain the theorem from Lemma 3.1. □

THEOREM 3.3. *There exist $g \in \mathcal{G}$ and $C \in (0, \infty)$ satisfying (3.1); equivalently, there exists $g \in \mathcal{G}$ such that*

$$(3.7) \qquad \limsup_{n \to \infty} \frac{a^n Y(-n)}{g(n)} = 1 \qquad Q\text{-a.s.}$$

*if and only if $K(x) \notin \mathcal{D}$.*

PROOF. We see from Theorem 3.2 that if $K(x) \in \mathcal{D}$, then there does not exist $g \in \mathcal{G}$ satisfying (3.1) with $C \in (0, \infty)$. Next, suppose that $K(x) \notin \mathcal{D}$. Then by Lemma 2.7, $\bar{\rho}_Y(x) \notin \mathcal{D}$. Since the support of $\nu_W$ is unbounded, we see from (2.18) and (2.24) that the support of $\rho_Y$ is also unbounded. Thus there is a sequence $y_n \uparrow \infty$ as $n \to \infty$ satisfying $y_{n+1} \geq 2y_n$ and

$$2^{-n} \bar{\rho}_Y(y_n) \geq \bar{\rho}_Y(2y_n) > 0.$$

Let $a_0 \in (1, a)$. We can take a strictly increasing sequence $\{x_n\}_{n=0}^{\infty}$ in such a way that $x_0 = 0$, $x_n \uparrow \infty$ as $n \to \infty$, $g(x) = y_n$ on $[x_{2n}, x_{2n+1})$ with $1 \leq \bar{\rho}_Y(y_n)(x_{2n+1} - x_{2n} - 2) \leq 2$; further set $g(x) = b_n a_0^x$ on $[x_{2n+1}, x_{2n+2})$ satisfying $b_n > 0$, $g(x_{2n+1}) = y_n$ and $g(x_{2n+2}) = y_{n+1}$. Then we have

$$\sum_{n=0}^{\infty} Q(Y(0) - Y(-1) > g(n)) \geq \sum_{n=0}^{\infty} \bar{\rho}_Y(y_n)(x_{2n+1} - x_{2n} - 2) = \infty.$$

Define $J_n := \{k \in \mathbb{Z}_+ : x_{2n+1} \leq k < x_{2n+2}\}$. Since we find from Lemma 2.10 that for some $c_1 > 0$

$$\sum_{n=0}^{\infty} \sum_{k \in J_n} \bar{\rho}_Y(b_n a_0^k) \leq c_1 \sum_{n=0}^{\infty} \int_{y_n}^{y_{n+1}} \bar{\rho}_Y(y) \frac{dy}{y} = c_1 \int_{y_0}^{\infty} \bar{\rho}_Y(y) \frac{dy}{y} < \infty,$$

we see that

$$\sum_{n=0}^{\infty} Q(Y(0) - Y(-1) > 2g(n))$$

$$\leq \sum_{n=0}^{\infty} \left( 2^{-n} \bar{\rho}_Y(y_n)(x_{2n+1} - x_{2n} + 1) + \sum_{k \in J_n} \bar{\rho}_Y(b_n a_0^k) \right) < \infty.$$

Thus we obtain (3.1) with $C \in (0, \infty)$ from Lemma 3.1 and Theorem 3.1. By replacing $g$ with $Cg$, we have (3.7). □

REMARK 3.1. Suppose that $K(x) \notin \mathcal{D}$. Define $\phi(t) := t^\alpha g(|\log t|)$ by using $g$ in the proof of the above theorem. Then obviously $\phi \in \Phi$ with increasing $g$.

PROPOSITION 3.1. *Let $C \in [0, \infty)$ and $g \in \mathcal{G}$. Suppose that $g(x)$ is increasing on $\mathbb{R}_+$ and the inverse function $g^{-1}(x)$ is quasi-submultiplicative on $(g(0), \infty)$. If*

$$(3.8) \qquad E(Wg^{-1}(\theta W)) \begin{cases} < \infty, & \text{for } 0 < \theta < C^{-1}, \\ = \infty, & \text{for } \theta > C^{-1}, \end{cases}$$

*then (3.1) holds.*

PROOF. Note that, for $\delta := \theta^{-1} > 0$,

$$\int_0^\infty Q(Y(0) > \delta g(x))\, dx = E_Q(g^{-1}(\theta Y(0))) = E(Wg^{-1}(\theta W)).$$

Suppose that (3.8) holds for $0 \leq C < \infty$. Then we obtain that

$$\sum_{n=0}^\infty Q(Y(0) > \delta g(n)) < \infty \qquad \text{for } \delta > C,$$

and hence, for every $\ell \geq 1$,

$$(3.9) \qquad \sum_{n=0}^\infty Q(Y(0) - Y(-\ell) > \delta g(n)) < \infty \qquad \text{for } \delta > C.$$

We see from the quasi-submultiplicativity of $g^{-1}(x)$ that, for any $\varepsilon > 0$,

$$\begin{aligned}
E_Q(g^{-1}(\theta Y(0))) \\
&\leq c_1 E_Q(g^{-1}((1+\varepsilon)\theta(Y(0) - Y(-\ell)))g^{-1}(c_2 \theta Y(-\ell))) \\
&= c_1 E_Q(g^{-1}((1+\varepsilon)\theta(Y(0) - Y(-\ell))))E_Q(g^{-1}(c_2 \theta Y(-\ell))) \\
&= c_1 E_Q(g^{-1}((1+\varepsilon)\theta(Y(0) - Y(-\ell))))E_Q(g^{-1}(c_2 a^{-\ell}\theta(Y(0)))).
\end{aligned}$$

Since, for sufficiently large $\ell$, $E_Q(g^{-1}(c_2 a^{-\ell}\theta Y(0))) = E(W(g^{-1}(c_2 a^{-\ell} \times \theta W))) < \infty$, we have

$$E_Q(g^{-1}((1+\varepsilon)\theta(Y(0) - Y(-\ell)))) = \infty \qquad \text{for } \theta > C^{-1} \text{ and any } \varepsilon > 0,$$

that is,

$$(3.10) \qquad \sum_{n=0}^\infty Q(Y(0) - Y(-\ell) > \delta g(n)) = \infty \qquad \text{for } 0 < \delta < C.$$

Thus we obtain (3.1) from (3.9) and (3.10) thanks to Theorem 3.1. $\square$

PROPOSITION 3.2. *Under the same assumption as Theorem 1.3, we have*

$$\limsup_{n\to\infty} \frac{a^n Y(-n)}{(\log n)^{(\gamma-1)/\gamma}} = \tau^{(1-\gamma)/\gamma} \qquad Q\text{-}a.s., \tag{3.11}$$

*where $\tau$ is a positive constant determined by (1.12).*

PROOF. Let $g(x) := (\log(e \vee x))^{(\gamma-1)/\gamma}$. Then $g^{-1}(x) = \exp(x^{\gamma/(\gamma-1)})$ for $x > 1$. Thus we see from Lemma 2.3 that (3.8) holds for $C = \tau^{(1-\gamma)/\gamma}$ and thereby Proposition 3.1 can be applied. □

PROPOSITION 3.3. *Under the same assumption as Theorem 1.4, we have*

$$\limsup_{n\to\infty} \frac{a^n Y(-n)}{\log n} = \sigma^{-1} \qquad Q\text{-}a.s., \tag{3.12}$$

*where $\sigma$ is a positive constant given by (1.14).*

PROOF. Let $g(x) := \log(e \vee x)$. Then $g^{-1}(x) = \exp x$ for $x > 1$. Thus we see from Lemma 2.4 that (3.8) holds for $C = \sigma^{-1}$ and thereby Proposition 3.1 can be applied. □

PROPOSITION 3.4. *Under the same assumption as Theorem 1.5, we have, Q-a.s.*

$$\limsup_{n\to\infty} \frac{a^n Y(-n)}{n^b} = \begin{cases} \infty, & \text{for } b < b_0, \\ 0, & \text{for } b > b_0. \end{cases} \tag{3.13}$$

*Moreover, if $E(N^{\theta_0}) < \infty$ (resp. $= \infty$), then*

$$\limsup_{n\to\infty} \frac{a^n Y(-n)}{n^{b_0}} = 0 \qquad (\text{resp.} = \infty) \ Q\text{-}a.s.$$

PROOF. Let $\ell = 1$ and $g(x) := \delta x^b$ with $\delta > 0$ and $b > 0$. Then $g^{-1}(x) = (x/\delta)^{1/b}$ for $x \geq 0$. Note that (3.3) holds if and only if $E_Q((Y(0) - Y(-1))^{1/b}) = \infty$. Thus we see from Lemma 2.9 that, for any $\delta > 0$, (3.3) holds for $b < b_0$ and (3.5) does for $b > b_0$. Therefore we obtain (3.13) from Lemma 3.1 for $b > 0$ and thereby also for $b \leq 0$. The second assertion can be proved in the same manner. □

PROPOSITION 3.5. *Under the same assumption as Theorem 1.6, we have*

$$\limsup_{n\to\infty} \frac{a^n Y(-n)}{R^{-1}(\log n)} = \xi_R^{-1} \qquad Q\text{-}a.s., \tag{3.14}$$

*where $\xi_R$ is a positive constant given by (1.18) and (1.18) is clearly equivalent to*

$$E(W \exp(R(\delta W))) \begin{cases} < \infty, & \text{for } 0 < \delta < \xi_R, \\ = \infty, & \text{for } \delta > \xi_R. \end{cases} \tag{3.15}$$

PROOF. Let $g(x) := R^{-1}(\log(e \vee x))$. Then $g^{-1}(x) = \exp(R(x))$ almost everywhere for $x > R^{-1}(1)$. It is obvious that $g^{-1}(x)$ satisfies quasi-subexponentiality. Thus we see from (3.15) that (3.8) holds for $C = \xi_R^{-1}$ and thereby Proposition 3.1 can be applied. $\square$

REMARK 3.2. Propositions 3.2 and 3.3 were proved by Liu [23]. However, each one-half of their proofs contained a serious gap because they depend on Theorems 1.3 and 1.4, respectively. See Remark 1.4. The first assertion of Proposition 3.4 was conjectured by Liu [21]. On the other hand, Proposition 3.5 was conjectured by Hawkes [17]. Many other limit theorems for the sequence $\{Y(-n)\}$ are found in [15, 23, 24, 33, 34, 38] and [46]. Further, limit theorems for another increasing shift self-similar additive random sequence associated with a Galton–Watson branching process are discussed in [45]. The study of normalizability type theorems such as Corollary 2.2 of [46] and Theorem 3.3 was motivated by a celebrated paper [35] for an increasing random walk.

**4. Proof of the main theorems.** In this section, we prove the main and additional theorems stated in Section 1 by using the results in Section 3. Lemmas 4.1 and 4.2 are crucial density theorems for the measure $\phi$-$H$. The first one is suggested by Proposition 3.2 of [21]. As was pointed out in general by Taylor [40], the second one is difficult to express without the exceptional set $\Delta$.

LEMMA 4.1. *Let $C \in [0, \infty)$ and $\phi(t) := t^\alpha g(|\log t|) \in \Phi$. Let $\Delta'$ be an arbitrary Borel set in $\partial \mathbf{T}$ with $\mu(\Delta') = 0$ a.s. If*

$$\limsup_{n \to \infty} \frac{a^n Y(-n)}{g(n)} \leq C \qquad Q\text{-a.s.,} \tag{4.1}$$

*then*

$$\phi\text{-}H(\partial \mathbf{T} \setminus \Delta') \geq C^{-1} W \qquad a.s. \text{ on } \{\partial \mathbf{T} \neq \varnothing\}, \tag{4.2}$$

*with the understanding that $0/0 = 0$.*

PROOF. Suppose that (4.1) holds with $0 \leq C < \infty$. For any $\varepsilon > 0$, there are compact set $K = K(\omega) \subset \partial \mathbf{T} \setminus \Delta'$ and a positive integer $n_0 = n_0(\omega)$ almost surely such that $\mu(K) \geq W - \varepsilon$ and

$$\mu(B_{\mathbf{i}|n}) \leq (C + \varepsilon) \phi(|B_{\mathbf{i}|n}|) \qquad \text{for every } \mathbf{i} \in K \text{ and every } n \geq n_0.$$

Thus we see that almost surely

(4.3) $\mu(B_{\mathbf{i}|n} \cap K) \leq (C + \varepsilon) \phi(|B_{\mathbf{i}|n}|)$ for every $\mathbf{i} \in \mathbf{I}$ and every $n \geq n_0$.

Let $\{S_j\}_{j=0}^{\infty}$ be an arbitrary cover of the set $K$ with $S_j$ being balls satisfying $|S_j| \leq e^{-n_0}$. Then, since we see from (4.3) that almost surely $\mu(S_j \cap K) \leq (C+\varepsilon)\phi(|S_j|)$,

$$W - \varepsilon \leq \mu(K) \leq \mu\left(\bigcup_{j=0}^{\infty}(S_j \cap K)\right)$$

$$\leq \sum_{j=0}^{\infty} \mu(S_j \cap K) \leq (C+\varepsilon)\sum_{j=0}^{\infty} \phi(|S_j|).$$

Thereby we obtain from the definition of Hausdorff measure that almost surely

$$\phi\text{-}H(\partial\mathbf{T} \setminus \Delta') \geq \phi\text{-}H(K) \geq \frac{W-\varepsilon}{C+\varepsilon}.$$

Letting $\varepsilon \downarrow 0$, we have (4.2). $\square$

LEMMA 4.2. *Let $C \in (0, \infty]$ and $\phi(t) := t^\alpha g(|\log t|) \in \Phi$. If*

$$\limsup_{n \to \infty} \frac{a^n Y(-n)}{g(n)} \geq C \qquad Q\text{-a.s.,} \tag{4.4}$$

*then we have*

$$\phi\text{-}H(\partial\mathbf{T} \setminus \Delta) \leq C^{-1}W \qquad \text{with } \mu(\Delta) = 0 \text{ a.s.} \tag{4.5}$$

*Thus, if $\phi$-$H$ is absolutely continuous with respect to $\mu$, then*

$$\phi\text{-}H(\partial\mathbf{T}) \leq C^{-1}W \qquad \text{a.s.} \tag{4.6}$$

PROOF. Suppose that (4.4) holds for $0 < C \leq \infty$. Define the set $\Delta(\delta)$ in $\partial\mathbf{T}$ for $\delta \geq 0$ as

$$\Delta(\delta) := \left\{\mathbf{i} \in \partial\mathbf{T} : \limsup_{n\to\infty} \frac{W_{\mathbf{i}|n}}{g(n)} \leq \delta\right\}.$$

Note that $\Delta = \Delta(0)$. We prove that

$$\phi\text{-}H(\Delta(C-) \setminus \Delta) = 0 \qquad \text{a.s.,} \tag{4.7}$$

with the understanding that $\Delta(\infty-) = \lim_{\delta\uparrow\infty}\Delta(\delta)$. We see from (4.4) that

$$\mu(\Delta) = \mu(\Delta(C-)) = 0 \qquad \text{a.s.} \tag{4.8}$$

Let $0 < \delta < C$ and $k \in \mathbb{N}$. Then define $m = m(\mathbf{i}) \geq k$ for $\mathbf{i} \in (\Delta(\delta))^c \cap \partial\mathbf{T}$ and $\Gamma_k$ as

$$m := \inf\{n \geq k : W_{\mathbf{i}|n} > \delta g(n)\} \quad \text{and} \quad \Gamma_k := \{\mathbf{i}|m(\mathbf{i}) : \mathbf{i} \in (\Delta(\delta))^c \cap \partial\mathbf{T}\}.$$

Define $\Gamma'_k$ and $\widetilde{\Gamma}_k$ as

$$\Gamma'_k := \{\mathbf{i}|k : \mathbf{i} \in \Delta(\delta)\} \quad \text{and} \quad \widetilde{\Gamma}_k := \Gamma_k \cup \Gamma'_k.$$

Then $\widetilde{\Gamma}_k$ is a cutset for $\partial \mathbf{T}$. Define $\Gamma_k(\mathbf{j})$ and $\widetilde{\Gamma}_k(\mathbf{j})$ for $\mathbf{j} \in \mathbf{U}$ with $|\mathbf{j}| = k$ as

$$\Gamma_k(\mathbf{j}) := \{\mathbf{i} : \mathbf{j} \leq \mathbf{i} \in \Gamma_k\} \quad \text{and} \quad \widetilde{\Gamma}_k(\mathbf{j}) := \{\mathbf{i} : \mathbf{j} \leq \mathbf{i} \in \widetilde{\Gamma}_k\}.$$

Then $\widetilde{\Gamma}_k(\mathbf{j})$ is a cutset for $\partial \mathbf{T} \cap B_{\mathbf{j}}$. Thus we have by (2.5)

$$\phi\text{-}H((\Delta(\delta))^c \cap B_{\mathbf{j}}) \leq \sum_{\mathbf{i} \in \Gamma_k(\mathbf{j})} a^{-|\mathbf{i}|} g(|\mathbf{i}|)$$

$$\leq \delta^{-1} \sum_{\mathbf{i} \in \widetilde{\Gamma}_k(\mathbf{j})} a^{-|\mathbf{i}|} W_{\mathbf{i}} = \delta^{-1} a^{-k} W_{\mathbf{j}} = \delta^{-1} \mu(B_{\mathbf{j}}) \qquad \text{a.s.}$$

Hence we have for a Borel set $A \subset \partial \mathbf{T}$

(4.9) $\qquad \phi\text{-}H((\Delta(\delta))^c \cap A) \leq \delta^{-1} \mu(A) \qquad$ a.s.

Setting $A = \Delta(C-)$ and then letting $\delta \downarrow 0$ in (4.9), we obtain (4.7) from (4.8). Letting $\delta \uparrow C$ and setting $A = \partial \mathbf{T}$ in (4.9), we have

$$\phi\text{-}H((\Delta(C-))^c \cap \partial \mathbf{T}) \leq C^{-1} \mu(\partial \mathbf{T}) = C^{-1} W \qquad \text{a.s.}$$

Thus we conclude (4.5) from (4.7) and (4.8), and thereby also see (4.6). $\square$

PROPOSITION 4.1. *Let $\phi(t) := t^\alpha g(|\log t|) \in \Phi$ and let $\Delta'$ be an arbitrary Borel set in $\partial \mathbf{T}$ with $\mu(\Delta') = 0$ a.s. Suppose that (3.1) holds with $C \in [0, \infty]$.*

(i) *If $0 < C < \infty$, then $0 < \phi\text{-}H(\partial \mathbf{T} \setminus \Delta) = C^{-1} W < \infty$ with $\mu(\Delta) = 0$ a.s. on $\{\partial \mathbf{T} \neq \varnothing\}$. Moreover, for $A \in \mathcal{B}(\mathbf{I})$ satisfying $A \subset \partial \mathbf{T} \setminus \Delta$, $\phi\text{-}H(A) = C^{-1} \mu(A)$ a.s.*

(ii) *If $C = 0$, then $\phi\text{-}H(\partial \mathbf{T}) = \phi\text{-}H(\partial \mathbf{T} \setminus \Delta') = \infty$ a.s. on $\{\partial \mathbf{T} \neq \varnothing\}$.*

(iii) *If $C = \infty$, then $\phi\text{-}H(\partial \mathbf{T} \setminus \Delta) = 0$ with $\mu(\Delta) = 0$ a.s.*

PROOF. Suppose that (3.1) holds with $C \in [0, \infty]$. First we prove assertion (i). If $0 < C < \infty$, then we see from Lemmas 4.1 and 4.2 that

$$0 < \phi\text{-}H(\partial \mathbf{T} \setminus \Delta) = C^{-1} W < \infty \qquad \text{with } \mu(\Delta) = 0 \text{ a.s. on } \{\partial \mathbf{T} \neq \varnothing\}.$$

The second assertion follows from the fact that we can find in the same way as the proof of (4.9) that, for all $\mathbf{i} \in \mathbf{U}$,

$$\phi\text{-}H(\partial \mathbf{T} \cap B_{\mathbf{i}} \setminus \Delta) = C^{-1} \mu(B_{\mathbf{i}}).$$

The proofs of assertions (ii) and (iii) are as follows. If $C = 0$, then we find from Lemma 4.1 that

$$\phi\text{-}H(\partial \mathbf{T}) = \phi\text{-}H(\partial \mathbf{T} \setminus \Delta') = \infty \qquad \text{a.s. on } \{\partial \mathbf{T} \neq \varnothing\}.$$

If $C = \infty$, then we obtain from Lemma 4.2 that

$$\phi\text{-}H(\partial \mathbf{T} \setminus \Delta) = 0 \qquad \text{with } \mu(\Delta) = 0 \text{ a.s.} \qquad \square$$

LEMMA 4.3. *Let $\phi(t) := t^\alpha g(|\log t|) \in \Phi$. If the condition $(G_\Delta)$ for $g$ is satisfied, then $\phi\text{-}H(\Delta) = 0$ a.s.*

PROOF. Thanks to $(G_\Delta)$, we can define a positive integer $n(k) > k$ such that, for some $\delta_0 > 0$,

$$(4.10) \quad \limsup_{k \to \infty} \left( \sum_{n=k}^{n(k)-1} Q(Y(0) - Y(-1) > \delta_0 g(n)) - \log g(n(k)) \right) = \infty.$$

Define the sets $A_k$ and $A'_k$ for $k \geq 1$ as

$$A_k := \{\mathbf{i} \in \partial \mathbf{T} : W_{\mathbf{i}|n} - a^{-1} W_{\mathbf{i}|(n+1)} \leq \delta_0 g(n) \text{ for all } n \geq k\}$$

and

$$A'_k := \{\mathbf{i} \in \partial \mathbf{T} : W_{\mathbf{i}|n} - a^{-1} W_{\mathbf{i}|(n+1)} \leq \delta_0 g(n) \text{ for } k \leq n \leq n(k) - 1\}.$$

Then we define $\Gamma_k$ and $\Gamma'_k$ as

$$\Gamma_k := \{\mathbf{i}|n(k) : \mathbf{i} \in A_k\} \quad \text{and} \quad \Gamma'_k := \{\mathbf{i}|n(k) : \mathbf{i} \in A'_k\}.$$

Then $\Gamma_k$ is a cutset for $A_k$. We observe that

$$(4.11) \quad E\left( \sum_{\mathbf{i} \in \Gamma'_k} W_{\mathbf{i}} \right) = E\left( \sum_{\mathbf{i} \in \Gamma'_k} 1 \right).$$

Let $\mathbf{m} := (m_0, m_1, \ldots, m_{n(k)-1}) \in \mathbb{N}^{n(k)}$ and define events $E_k := E_k(\mathbf{i}, \mathbf{m})$ and $H_k := H_k(\mathbf{i}, \mathbf{m})$ for $\mathbf{i} \in \mathbb{Z}_+^{n(k)}$ and $\mathbf{m} \in \mathbb{N}^{n(k)}$ as

$$E_k := \{\omega \in \Omega : N_{\mathbf{i}|n} = m_n \text{ for } n = 0, 1, \ldots, n(k) - 1\}$$

and

$$H_k := \left\{ \omega \in \Omega : \sum_{i=0}^{*(m_n - 1)} a^{-1} W_{(\mathbf{i}|n)*i} \leq \delta_0 g(n) \text{ for } k \leq n \leq n(k) - 1 \right\},$$

where the symbol $\sum_{i=0}^{*(m_n-1)}$ denotes the sum over $i$ from $0$ to $m_n - 1$ except for $i$ satisfying $W_{(\mathbf{i}|n)*i} = W_{\mathbf{i}|(n+1)}$. Define $G_k := \{\mathbf{i} \in \mathbb{Z}_+^{n(k)} : 0 \leq i_{n+1} \leq m_n - 1 \text{ for } 0 \leq n \leq n(k) - 1\}$. By using (2.4) in the second equality, we have

$$E\left( \sum_{\mathbf{i} \in \Gamma'_k} W_{\mathbf{i}} \right) = \sum_{\mathbf{m} \in \mathbb{N}^{n(k)}} \sum_{\mathbf{i} \in \mathbb{Z}_+^{n(k)}} E(\mathbb{1}_{E_k}(\omega) \mathbb{1}_{\Gamma'_k}(\mathbf{i}) W_{\mathbf{i}})$$

$$= \sum_{\mathbf{m} \in \mathbb{N}^{n(k)}} \sum_{\mathbf{i} \in G_k} E(\mathbb{1}_{H_k \cap E_k}(\omega) W_{\mathbf{i}}).$$

Since $H_k$, $E_k$ and $W_\mathbf{i}$ for $\mathbf{i} \in G_k$ are independent, we see from $E(W_\mathbf{i}) = 1$ for $\mathbf{i} \in G_k$ that

$$E\left(\sum_{\mathbf{i}\in\Gamma'_k} W_\mathbf{i}\right) = \sum_{\mathbf{m}\in\mathbb{N}^{n(k)}} \sum_{\mathbf{i}\in G_k} E(\mathbb{1}_{H_k \cap E_k}(\omega))E(W_\mathbf{i})$$

$$= \sum_{\mathbf{m}\in\mathbb{N}^{n(k)}} \sum_{\mathbf{i}\in\mathbb{Z}_+^{n(k)}} E(\mathbb{1}_{E_k}(\omega)\mathbb{1}_{\Gamma'_k}(\mathbf{i}))$$

$$= E\left(\sum_{\mathbf{i}\in\Gamma'_k} 1\right).$$

Thus we have proved (4.11). By using (4.11), the expectation of $\phi$-$H(A_k)$ is estimated as

$E(\phi\text{-}H(A_k))$

$$\leq E\left(\sum_{\mathbf{i}\in\Gamma_k} \phi(|B_\mathbf{i}|)\right) = E\left(\sum_{\mathbf{i}\in\Gamma_k} a^{-n(k)}g(n(k))\right)$$

$$\leq E\left(\sum_{\mathbf{i}\in\Gamma'_k} a^{-n(k)}g(n(k))\right) = g(n(k))E\left(\sum_{\mathbf{i}\in\Gamma'_k} a^{-n(k)}W_\mathbf{i}\right)$$

$$= g(n(k))Q(Y(-n) - Y(-n-1) \leq a^{-n}\delta_0 g(n) \text{ for } k \leq n \leq n(k)-1)$$

$$= g(n(k)) \prod_{n=k}^{n(k)-1} Q(Y(-n) - Y(-n-1) \leq a^{-n}\delta_0 g(n))$$

$$\leq \exp\left(-\sum_{n=k}^{n(k)-1} Q(Y(0) - Y(-1) > \delta_0 g(n)) + \log g(n(k))\right),$$

where we used the additivity and shift self-similarity of the sequence $\{Y(n)\}$ in the last equality and inequality, respectively. Note that $\Delta \subset \lim_{k\to\infty} A_k = \bigcup_{k=1}^{\infty} A_k$. Thus we see from (4.10) that

$$E(\phi\text{-}H(\Delta)) \leq \liminf_{k\to\infty} E(\phi\text{-}H(A_k)) = 0.$$

Therefore we have $\phi$-$H(\Delta) = 0$ a.s. $\square$

LEMMA 4.4. *Let $C \in (0,\infty)$ and $\phi(t) := t^\alpha g(|\log t|) \in \Phi$ with increasing $g \in \mathcal{G}$. Suppose that (3.1) holds. Then there exists $\phi^*(t) := t^\alpha g^*(|\log t|) \in \Phi$ such that $g^*(x) \leq g(x)$ for $x \geq 0$, $g^*$ satisfies the condition ($G_\Delta$) and*

(4.12) $$\limsup_{n\to\infty} \frac{a^n Y(-n)}{g^*(n)} = C \qquad Q\text{-a.s.}$$

PROOF. Without harming (3.1), we can and do assume from Theorem 3.1 that $g \in \mathcal{G}$ is left-continuous adding to be increasing. However, $g^*(x) \in \mathcal{G}$ defined below is not always increasing and left-continuous. Let $\{x_n\}_{n=0}^{\infty}$ be an increasing sequence satisfying $x_0 = 0$, $\lim_{n \to \infty} x_n = \infty$, and $x_{2n+1} \in \mathbb{Z}_+$ for $n \in \mathbb{Z}_+$. Let $I_n := \{k \in \mathbb{Z}_+ : x_{2n} \le k < x_{2n+1}\}$ and $J_n := \{k \in \mathbb{Z}_+ : x_{2n+1} \le k \le x_{2n+2}\}$ for $n \in \mathbb{Z}_+$. Define $G(n)$ for $n \in \mathbb{Z}_+$ as

$$G(n) := \sum_{k \in I_n} Q(Y(0) - Y(-1) > \delta_0 g(k))$$

with sufficiently small $\delta_0 > 0$ satisfying $\sum_{n=0}^{\infty} Q(Y(0) - Y(-1) > \delta_0 g(n)) = \infty$ owing to (ii) of Lemma 3.1. Let $a_0 := \limsup_{n \to \infty} g(n+1)/g(n)$ and put $a_1$ as $a_1 \in (a_0, a)$. We choose the sequence $\{x_n\}_{n=0}^{\infty}$ in such a way that $x_{2n} < x_{2n+1} \le x_{2n+2}$ and $e^{G(n)} \ge g(x_{2n})$ for $n \in \mathbb{Z}_+$. Up to this step, there is freedom of the choice of $x_{2n+2}$ except for $x_{2n+1} \le x_{2n+2}$ for $n \in \mathbb{Z}_+$. Next set $g^* \in \mathcal{G}$ as follows. For $n \in \mathbb{Z}_+$, $g^*(x) := g(x)$ on $x_{2n} \le x < x_{2n+1}$, $g^*(x_{2n+1}) := e^{G(n)} \wedge g(x_{2n+1})$, and $g^*(x) := b_n a_1^x$ on $x_{2n+1} \le x \le x_{2n+2}$ with some $b_n > 0$. Further, the equality $g^*(x_{2n+2}) := g(x_{2n+2})$ and the inequality $g^*(x) \le g(x)$ on $x_{2n+1} \le x \le x_{2n+2}$ are possible by defining $x_{2n+2} := \sup\{x \ge x_{2n+1} : b_n a_1^x \le g(x)\}$ since $g$ is left-continuous and the set $\{x \ge x_{2n+1} : b_n a_1^x \le g(x)\}$ is nonempty and bounded by virtue of $a_1 > a_0 \ge 1$. Note that $g^*(x_{2n}) \le g^*(x_{2n+1})$ for $n \in \mathbb{Z}_+$ and that the equality $x_{2n+1} = x_{2n+2}$ can hold in case $e^{G(n)} \ge g(x_{2n+1})$. Then we see that

$$\limsup_{n \to \infty} \left( \sum_{k=0}^{x_{(2n+1)}-1} Q(Y(0) - Y(-1) > \delta_0 g^*(k)) - \log g^*(x_{2n+1}) \right)$$

$$\ge \sum_{n=0}^{\infty} Q(Y(0) - Y(-1) > \delta_0 g(n)) = \infty.$$

Thus $(G_\Delta)$ holds for $g^* \in \mathcal{G}$. Moreover, we find that, for $\delta > C$,

$$\sum_{n=0}^{\infty} \sum_{k \in I_n} Q(Y(0) - Y(-\ell) > \delta g^*(k))$$

$$= \sum_{n=0}^{\infty} \sum_{k \in I_n} Q(Y(0) - Y(-\ell) > \delta g(k)) < \infty$$

and from Lemma 2.10 that with some $c_1 > 0$

$$\sum_{n=0}^{\infty} \sum_{k \in J_n} Q(Y(0) - Y(-\ell) > \delta g^*(k))$$

$$= \sum_{n=0}^{\infty} \sum_{k \in J_n} Q(Y(0) - Y(-\ell) > \delta b_n a_1^k)$$

$$\leq c_1 \int_{\delta g(0)}^{\infty} Q(Y(0) - Y(-\ell) > y) \frac{dy}{y} < \infty.$$

Thus we obtain that, for any $\delta > C$ and any $\ell \geq 1$,

$$\sum_{n=0}^{\infty} Q(Y(0) - Y(-\ell) > \delta g^*(n))$$

$$\leq \sum_{n=0}^{\infty} \left( \sum_{k \in I_n} + \sum_{k \in J_n} \right) Q(Y(0) - Y(-\ell) > \delta g^*(k)) < \infty.$$

On the other hand, for any $\delta \in (0, C)$ and some $\ell \geq 1$,

$$\sum_{n=0}^{\infty} Q(Y(0) - Y(-\ell) > \delta g^*(n))$$

$$\geq \sum_{n=0}^{\infty} Q(Y(0) - Y(-\ell) > \delta g(n)) = \infty.$$

Therefore we have established (4.12) by Theorem 3.1. □

PROOF OF THEOREM 1.1. Suppose that $K(x) \notin \mathcal{D}$. Then we see from Theorem 3.3, Remark 3.1 and Lemma 4.4 that there is $\phi^* \in \Phi$ satisfying the conditions $(G_\Delta)$, (3.1) and (3.2) with $g = g^*$ and $C \in (0, \infty)$. Replacing $\phi$ by $\phi^*$, we obtain (1.8) from Proposition 4.1 and Lemma 4.3 with $C_\phi = C^{-1} \in (0, \infty)$. The second assertion follows from Proposition 4.1. □

PROOF OF THEOREM 1.2. Suppose that $K(x) \in \mathcal{D}$ and let $\phi \in \Phi$. We first prove assertion (i). If $\sum_{n=0}^{\infty} K(g(n)) < \infty$, then we see from Theorem 3.2 that (3.1) holds with $C = 0$. Thus we obtain assertion (i) from (ii) of Proposition 4.1. Next we prove assertion (ii). If $\sum_{n=0}^{\infty} K(g(n)) = \infty$, then we find from Theorem 3.2 that (3.1) holds with $C = \infty$. Thus we obtain assertion (ii) from (iii) of Proposition 4.1. Lastly we prove assertion (iii). Suppose that $\sum_{n=0}^{\infty} K(g(n)) = \infty$ and $\limsup_{\delta \to 0+} \limsup_{n \to \infty} \sum_{k=0}^{n} K(\delta g(k))/\log(e \vee g(n)) = \infty$. Then we observe that $(G_\Delta)$ holds. Note from Lemma 2.7 that

$$Q(Y(0) - Y(-1) > x) \asymp K(x) \quad \text{as } x \to \infty.$$

Let $\delta_0 > 0$. In the case where $M_0 := \liminf_{n \to \infty} g(n) < \infty$, we find from $\sum_{n=0}^{\infty} K(g(n)) = \infty$ that

$$\limsup_{n \to \infty} \left( \sum_{k=0}^{n-1} Q(Y(0) - Y(-1) > \delta_0 g(k)) - \log g(n) \right)$$

$$= \sum_{n=0}^{\infty} Q(Y(0) - Y(-1) > \delta_0 g(n)) - \log M_0 = \infty.$$

In the case where $\lim_{n\to\infty} g(n) = \infty$, we see that

$$\limsup_{\delta_0 \to 0+} \limsup_{n\to\infty} \frac{\sum_{k=0}^{n-1} Q(Y(0) - Y(-1) > \delta_0 g(k))}{\log g(n)}$$

$$= \limsup_{\delta_0 \to 0+} \limsup_{n\to\infty} \frac{\sum_{k=0}^{n-1} K(\delta_0 g(k))}{\log g(n)} = \infty.$$

Thus $(G_\Delta)$ holds. It follows from Lemma 4.3 that $\phi\text{-}H(\Delta) = 0$ a.s. and thereby combining with assertion (ii), we conclude that

$$\phi\text{-}H(\partial \mathbf{T}) = \phi\text{-}H(\partial \mathbf{T} \setminus \Delta) = 0 \qquad \text{a.s.} \qquad \square$$

PROOF OF COROLLARY 1.1. The corollary is clear from Theorems 1.1 and 1.2 and Proposition 4.1. In particular, the last assertion of (ii) is obvious from Lemma 4.3. $\square$

The following lemma is a special case of Lemma 3.2 of [21].

LEMMA 4.5. *Suppose that $h(x)$ is a nonnegative decreasing function on $\mathbb{R}_+$ with $\int_1^\infty h(x)\,dx = \infty$. Then we have, for all $\delta > 0$ and for all $\varepsilon_1 \in (0, \delta)$,*

$$(4.13) \qquad \limsup_{n\to\infty} \left( \int_1^{n^{1/(1+\delta)}} h(x) x^\delta \, dx - n^{\varepsilon_1/(1+\delta)} \right) = \infty.$$

The proofs of Theorems 1.3, 1.4 and 1.6 are due to the following proposition.

PROPOSITION 4.2. *Let $C \in (0, \infty)$ and $\phi(t) := t^\alpha g(|\log t|) \in \Phi$. Suppose that $g(x)$ is increasing on $\mathbb{R}_+$ and the inverse function $g^{-1}(x)$ is quasi-submultiplicative on $(g(0), \infty)$ and that there is $b \in (a^{-1}, 1)$ such that*

$$(4.14) \qquad \limsup_{x\to\infty} \frac{\log g^{-1}(bx)}{\log g^{-1}(x)} < 1.$$

*If (3.8) is satisfied, then (1.8) and (1.10) hold with $C_\phi = C^{-1}$.*

PROOF. Thanks to Propositions 3.1 and 4.1 and Lemma 4.3, it is enough to prove that $(G_\Delta)$ holds for $g$. Choosing $c_1 > 0$ satisfying $a^{-1} c_1^{-1} < C^{-1} < b c_1^{-1}$, we have by (3.8)

$$(4.15) \qquad \int_0^\infty \bar{\eta}_Y(ac_1 g(x))dx = E(W g^{-1}(a^{-1} c_1^{-1} W)) < \infty$$

and

$$(4.16) \qquad \int_0^\infty \bar{\eta}_Y(b^{-1} c_1 g(x))\,dx = E(W g^{-1}(b c_1^{-1} W)) = \infty.$$

Note from (4.14) that there are $\delta_1 \in (0,1)$, $M_1 > 0$ and positive constants $c_2$ and $c_3$ such that

(4.17)
$$g^{-1}(bx) \leq (g^{-1}(x))^{\delta_1} \quad \text{and}$$
$$g(x) \leq c_2 (\log x)^{c_3} \qquad \text{for } x > M_1.$$

Thus we obtain from (2.18), (4.15) and (4.17) with $\delta_0 := \varepsilon c_1$ and $\delta := \delta_1^{-1} - 1$ that, for sufficiently large $n$,

$$\sum_{k=0}^{n-1} Q(Y(0) - Y(-1) > \delta_0 g(k))$$
$$\geq \int_1^n \bar{\rho}_Y(\varepsilon c_1 g(x)) \, dx$$
$$\geq 2^{-1} \left( \int_1^n \bar{\eta}_Y(c_1 g(x)) \, dx - \int_1^n \bar{\eta}_Y(ac_1 g(x)) \, dx \right)$$
$$\geq c_4 \int_1^n Q(g^{-1}(c_1^{-1} Y(0)) > x) \, dx$$
$$\geq c_4 \int_1^n Q(g^{-1}(bc_1^{-1} Y(0)) > x^{\delta_1}) \, dx$$
$$= \delta_1^{-1} c_4 \int_1^{n^{\delta_1}} Q(g^{-1}(bc_1^{-1} Y(0)) > y) y^{\delta} \, dy$$

with some positive constant $c_4$. Set $h(x) := Q(g^{-1}(bc_1^{-1} Y(0)) > x)$. Then we see from (4.16) that $\int_1^\infty h(x) \, dx = \infty$. It follows from (4.17) and Lemma 4.5 that

$$\limsup_{n\to\infty} \left( \sum_{k=0}^{n-1} Q(Y(0) - Y(-1) > \delta_0 g(k)) - \log g(n) \right)$$
$$\geq \limsup_{n\to\infty} \left( \delta_1^{-1} c_4 \int_1^{n^{\delta_1}} h(x) x^\delta \, dx - \log g(n) \right)$$
$$\geq \infty + \liminf_{n\to\infty} (\delta_1^{-1} c_4 n^{\varepsilon_1 \delta_1} - \log(c_2 (\log n)^{c_3})) = \infty.$$

Thus we have established $(G_\Delta)$ holds for $g$. □

PROOF OF THEOREM 1.3. Let $g(x) := (\log(e \vee x))^{(\gamma-1)/\gamma}$. Then $g^{-1}(x) = \exp(x^{\gamma/(\gamma-1)})$ for $x > 1$. Thus Proposition 4.2 can be applied with $C = \tau^{(1-\gamma)/\gamma}$. □

PROOF OF THEOREM 1.4. Let $g(x) := \log(e \vee x)$. Then $g^{-1}(x) = \exp x$ for $x > 1$. Thus Proposition 4.2 can be applied with $C = \sigma^{-1}$. □

PROOF OF THEOREM 1.5. Let $g(x) := x^b$ on $\mathbb{R}_+$. If $b > b_0$, then $\psi_b\text{-}H(\partial \mathbf{T}) = \infty$ by Propositions 3.4 and 4.1. If $b < b_0$, then choose $\delta_1 \in (0,1)$ satisfying $b\delta_1^{-1} < b_0$ and set $\delta := \delta_1^{-1} - 1$ and $h(x) := Q(Y(0) - Y(-1) > x^{b\delta_1^{-1}})$. Note from Lemma 2.9 that $\int_1^\infty h(x)\,dx = \infty$. Thus we obtain from Lemma 4.5 that, for $b < b_0$ and $0 < \varepsilon_1 < \delta$,

$$\limsup_{n\to\infty}\left(\sum_{k=0}^{n-1} Q(Y(0) - Y(-1) > g(k)) - \log g(n)\right)$$
$$\geq \limsup_{n\to\infty}\left(\int_1^n Q(Y(0) - Y(-1) > y^b)\,dy - \log n^b\right)$$
$$= \limsup_{n\to\infty}\left(\delta_1^{-1}\int_1^{n^{\delta_1}} h(x)x^\delta\,dx - b\log n\right)$$
$$\geq \infty + \liminf_{n\to\infty}(\delta_1^{-1} n^{\varepsilon_1\delta_1} - b\log n) = \infty.$$

Hence $(G_\Delta)$ holds and thereby $\psi_b\text{-}H(\partial \mathbf{T}) = \psi_b\text{-}H(\partial \mathbf{T}\setminus\Delta) = 0$ a.s. by Propositions 3.4 and 4.1 and Lemma 4.3. The second assertion is obvious from Propositions 3.4 and 4.1. □

PROOF OF THEOREM 1.6. Let $g(x) := R^{-1}(\log(e \vee x))$. Then $g^{-1}(x) = \exp(R(x))$ almost everywhere for $x > R^{-1}(1)$. It is obvious that $g^{-1}(x)$ satisfies the assumptions of Proposition 4.2. Thus Proposition 4.2 can be applied with $C = \xi_R^{-1}$. □

CONCLUDING REMARKS. Hawkes [17] proposed an outline for the resolution of the problem of determining the exact Hausdorff measure on the boundary of a Galton–Watson tree. The first step is to study "limsup" type limit theorems for the sequence $\{Y(n)\}$. The second step is to apply those limit theorems to determine the exact Hausdorff measure. This paper resolved the first step of his outline, but did not completely resolve the second step. Thus it is still unanswered whether there exists an exact Hausdorff measure which is not absolutely continuous with respect to the branching measure. The point is whether $\phi\text{-}H(\Delta) = 0$ a.s. for any exact Hausdorff measure $\phi\text{-}H$. However, the exceptional set $\Delta$ is so difficult to manage that the final goal might be beyond our way of approach. We end this article by posing the following problem which is an extension of Theorem 1.6 (Hawkes's conjecture).

PROBLEM. How is the exact Hausdorff measure explicitly given in the case where the distribution of $W$ is subexponential or O-subexponential?

This problem is deeply connected with the open problems in the Appendix of [45]. For the definitions of subexponentiality and O-subexponentiality, see [39].

**Acknowledgment.** The author is grateful to K. Sato for helpful advice and comments.

CENTER FOR MATHEMATICAL SCIENCES
UNIVERSITY OF AIZU
AIZU-WAKAMATSU 965-8580
JAPAN
E-MAIL: t-watanb@u-aizu.ac.jp